\documentclass[journal]{IEEEtran}
\usepackage{amsmath,amsfonts,bm}
\usepackage{array}
\usepackage{listings}
\usepackage{adjustbox}
\usepackage[group-digits=true,group-minimum-digits=4]{siunitx}
\usepackage[caption=false,font=normalsize,labelfont=sf,textfont=sf]{subfig}
\usepackage{tikz,pgfplots}

\usepackage{stfloats}
\usepackage{url}
\usepackage{verbatim}
\usepackage{graphicx}
\usepackage{acronym}

\usepackage{float}
\usepackage{comment}
\usepackage{soul}
\usepackage{color}
\usepackage{mathtools}
\usepackage{graphicx}
\usepackage{balance}
\usepackage{multicol}
\usepackage{longtable}
\usepackage{amssymb}
\usepackage{booktabs, makecell, tabularx}
\usepackage{longtable}
\usepackage[section]{placeins}
\usepackage{nomencl}
\usepackage{footmisc}
\usepackage{cite}
\usepackage{xcolor}
\pgfplotsset{compat=1.18}
\usepackage{setspace}
\usepackage{amsthm}
\usepackage{academicons}

\usepackage[acronym]{glossaries}
\usepackage{mathtools, cuted}
\usepackage{lipsum}
\usepackage{balance}
\usepackage{xpatch}
\usepackage{graphicx,color}
\usepackage{textcomp}
\usepackage[colorlinks,
            urlcolor=blue,
            linkcolor=blue,
            anchorcolor=blue,
            citecolor=blue
           ]{hyperref}
\theoremstyle{plain}
\newtheorem{theorem}{Theorem}    
\newtheorem{corollary}{Corollary}[theorem]

\makeatletter
\xpatchcmd{\proof}{\@addpunct{.}}{\@addpunct{:}}{}{}
\makeatother

\begin{document}

\title{SLIPT for Underwater IoT: System Modeling and Performance Analysis}

\author{Shunyuan Shang, Ziyuan Shi,~\IEEEmembership{Member,~IEEE}, and~Mohamed-Slim~Alouini,~\IEEEmembership{Fellow,~IEEE} 
\thanks{(\textit{Corresponding author: Ziyuan Shi})}
\thanks{S. Shang, Z. Shi and M.-S. Alouini are with the Computer, Electrical, and Mathematical Science and Engineering (CEMSE) Division, King Abdullah University of Science and Technology (KAUST), Thuwal, Makkah Province, Saudi Arabia (e-mail: shunyuan.shang@kaust.edu.sa;
ziyuan.shi@kaust.edu.sa; slim.alouini@kaust.edu.sa).}
\thanks{This work was supported by an ERIF/OSSARI Grant.}
}


\maketitle
\begin{abstract}
This paper presents a unified analytical framework for a two phase underwater wireless optical communication (UWOC) system that integrates Simultaneous Lightwave Information and Power Transfer (SLIPT) using a photovoltaic (PV) panel receiver. The proposed architecture enables self powered underwater sensor nodes by leveraging wide area and low cost PV panels for concurrent optical signal detection and energy harvesting. We develop a composite statistical channel that combines distance dependent absorption, turbulence induced fading characterized by the mixture Exponential Generalized Gamma (EGG )distribution, and beam misalignment due to pointing errors. Based on this model we derive closed form expressions for the probability density function, the cumulative distribution function, the outage probability (OP), the average bit error rate, the ergodic capacity, and the harvested power using Meijer G and Fox H functions. Overall, the paper introduces a practical analytical framework that provides clear guidance for design, optimization, and operation of SLIPT based UWOC systems.
\end{abstract}


\section{Introduction}

Underwater wireless optical communication (UWOC) has become an increasingly vital technology for achieving high speed, low latency, and secure data transmission in complex and dynamic underwater environments \cite{mohammed2024underwater,fang2023high}.In contrast to traditional acoustic methods, which face limited bandwidth, long propagation delays, and strong susceptibility to ambient noise and multipath distortion \cite{Stojanovic2007On}, UWOC leverages the optical domain to support data rates from several megabits to multiple gigabits per second \cite{Kaushal2016UWOC}.

{
Beyond link level benefits, underwater IoT networks are strategically important because they enable persistent sensing and actuation across large ocean spaces. Networks of fixed sensors, mobile gateways, and autonomous platforms such as AUVs provide continuous measurements of temperature, salinity, dissolved oxygen, pH, turbidity, and other biogeochemical indicators that support climate science, ecosystem management, and early warning for hazards such as earthquakes and tsunamis. In the industrial domain, they underpin inspection and condition monitoring of subsea pipelines, power cables, offshore wind farms, and oil and gas infrastructure, which improves safety and operational efficiency while reducing human risk. They also support aquaculture monitoring, resource mapping, maritime security, and collaborative scientific exploration, which motivates high throughput low latency links and energy autonomy at the node level \cite{Kao2017IoUT,Mohsan2023IoUT}.

Sustaining the energy supply of underwater IoT nodes remains a primary bottleneck \cite{alamu2023energy,melki2025auv}. Cabled power is impractical beyond short ranges because tethers add drag, restrict mobility, and reduce reliability in corrosive seawater. Periodic battery replacement is costly and often infeasible in remote or deep deployments. Inductive or acoustic power transfer experiences strong distance dependent loss. Far field radio frequency (RF) transfer is even less suitable in seawater because the high conductivity causes severe frequency dependent attenuation \cite{Theocharidis2025Underwater}. Within this context, laser power beaming in the blue green spectral window provides compact and directional energy delivery. The same optical flux can carry data while delivering energy through simultaneous lightwave information and power transfer (SLIPT) \cite{de2020toward}, which enables self powered operation of underwater sensors while preserving the high throughput and low latency of UWOC.

Three SLIPT strategies are commonly studied. Time switching divides the reception window into intervals for decoding and for harvesting. Power splitting allocates part of the received optical signal to a data path and part to a harvesting path in real time \cite{ma2019simultaneous}. Photovoltaic (PV) based reception replaces a conventional photodiode with a solar panel or applies spectral splitting so that one optical flux supports both information detection and energy harvesting \cite{maragliano2015demonstration,wang2015design,wang2014towards,kong2018underwater}. PV  receivers relax alignment through a large collection area and can scavenge energy from ambient or dedicated illumination. At the same time, junction capacitance and readout circuitry limit detection bandwidth, which requires interface design that balances information current and harvested direct current.
\begin{table*}[ht]
\label{tab:prior}
\centering
\caption{{ Comparison with representative prior work}}
\scriptsize
\setlength{\tabcolsep}{3.2pt}
\renewcommand{\arraystretch}{1.18}
{\begin{tabular}{p{1.3cm} p{4.8cm} p{3.5cm} p{5.0cm} p{1.5cm}}
\hline
Work  & System or method & SLIPT or Energy Harvesting (EH)   & Channel and misalignment & PV receiver \\
\hline
\cite{Wang2015JSAC}   & Indoor optical SLIPT with PV front end and circuit techniques for bandwidth extension & SLIPT at receiver & No underwater turbulence model. Misalignment not central & Yes \\
\cite{Shin2016OE}     & Terrestrial optical SLIPT with PV detection and interface design & SLIPT at receiver & No underwater statistics. Misalignment not central & Yes \\
\cite{Ammar2022CL}    & UWOC with power splitting and combined solar cell and SPAD front end & SLIPT & Energy aware operation. No closed form performance under composite turbulence and pointing errors & Yes (with SPAD) \\
\cite{Li2024AO}       & Parallel relay assisted UWOC & SLIPT & Energy efficiency optimization under aggregate fading. Misalignment is not central & No \\
\cite{Agarwal2024TETT}& Two hop UWOC and RF via buoy relay & SLIPT at relay & Hybrid optical and RF. Simplified treatment of misalignment & No \\
\cite{Ijeh2022JOCN}   & Vertical UWOC link & No SLIPT & Oceanic turbulence with pointing errors. Analytical outage evaluation & No \\
\cite{Zhang2024AO}    & Vertical UWOC with fixed gain amplify and forward relay & No SLIPT & Cascaded turbulence with generalized misalignment. Closed form outage and capacity with Meijer G and bivariate Fox H & No \\
\cite{Tong2022OL}     & High speed UWOC using series connected solar array detector & EH capable hardware & Experimental link. No unified statistical SLIPT analysis & Yes \\
This work             & Two phase UWOC with PV based SLIPT using power splitting & SLIPT & Beer Lambert attenuation with mixture EGG turbulence and Gaussian displacement pointing jitter. Closed form reliability expressions & Yes \\
\hline
\end{tabular}}
\end{table*}

Prior optical SLIPT studies with solar-panel receivers mainly target terrestrial and indoor links. They show that a solar panel can serve as both optical detector and energy harvester and that circuit techniques can extend detection bandwidth and boost extracted energy, yet they do not develop statistical models for underwater turbulence or misalignment \cite{Wang2015JSAC,Shin2016OE}.  Energy-aware UWOC SLIPT with power splitting and a combined solar-cell and single-photon avalanche diode front end has also been reported, although closed-form performance under composite turbulence and pointing errors is not derived \cite{Ammar2022CL}. Complementary UWOC works consider SLIPT without explicit PV modeling and emphasize energy efficiency in relay-assisted links \cite{Li2024AO}, while hybrid UWOC–RF systems execute SLIPT at an intermediate relay \cite{Agarwal2024TETT}. On the statistical side, vertical UWOC links with oceanic turbulence and pointing errors admit analytical outage evaluation \cite{Ijeh2022JOCN}, and relay-assisted vertical links with generalized misalignment yield compact Meijer G and bivariate Fox H expressions with closed-form outage and capacity \cite{Zhang2024AO}. Experiments further confirm the feasibility of PV-based detection through high-speed UWOC using series-connected solar arrays \cite{Tong2022OL}. We compare representative studies with our framework in Table~\ref{tab:prior}.}

Realizing this potential in practice requires a robust and physically consistent analytical framework that models underwater optical propagation and supports quantitative performance assessment. Early studies often relied on the Beer–Lambert law, which provides a simple exponential description of attenuation in homogeneous media \cite{cai2021analysis, ata2022absorption}. Although computationally efficient, this formulation does not capture the random intensity variations produced by aquatic turbulence. To better capture the random intensity variations produced by aquatic turbulence, the authors of \cite{Zedini2018EGG} introduced the Exponential Generalized Gamma (EGG) distribution from experimental measurements, which accurately characterizes channel statistics across different water qualities.

In this paper, we integrate a solar panel receiver into a two phase UWOC SLIPT architecture and place the focus on a tractable analytical formulation. We adopt a composite channel that merges Beer–Lambert attenuation with mixture EGG turbulence and Gaussian jitter pointing loss. This model yields exact closed form expressions for the probability density function and cumulative distribution function of the channel gain and the end to end signal to noise ratio. Building on these results, we derive outage probability (OP), average bit error rate, ergodic capacity, and average harvested power using Meijer G and Fox H functions. To validate the framework and explore sensitivity to environmental and design parameters, we conduct Monte Carlo simulations under varying turbidity levels, alignment conditions, wavelength choices, receiver apertures, and link ranges. The analysis provides design oriented guidance for link range planning, split factor selection, receiver aperture sizing, and device efficiency budgeting in practical deployments.

{ The primary contributions of this work are as follows:
\begin{itemize}
\item A unified analytical framework for SLIPT enabled UWOC that jointly incorporates distance dependent attenuation, a mixture EGG turbulence model, and Gaussian pointing jitter, together with photovoltaic reception and power splitting.
\item Exact closed form distributions and moments of the composite channel gain and the instantaneous SNR, enabling closed form expressions for OP, average bit error rate, ergodic capacity, and average harvested power.
\item A design oriented interpretation that maps analytical parameters to practical choices, including the power split factor, receiver aperture, jitter level, device efficiencies, wavelength, and link distance, with guidance on energy and reliability trade offs.
\item Monte Carlo validation that confirms the accuracy of the analysis and illustrates performance trends under representative underwater conditions.
\end{itemize}}

The remainder of the paper is organized as follows. Section~\ref{sec:system_model} describes the system architecture and signal model. Section~\ref{sec:channelmodel} develops the statistical channel model. Section~\ref{Perana} presents the performance analysis and closed form derivations. Section~\ref{sec:NUMERICAL ANALYSIS} reports simulation results and design insights. Section~\ref{sec:Conclusion} concludes the paper and outlines future research directions.

\section{System Model}
\label{sec:system_model}

The proposed UWOC system adopts a dual-mode PV energy-management architecture. During normal operation, a lithium-ion battery supplies sensing, processing, and housekeeping loads, which is consistent with energy-harvesting sensor nodes that use a battery to sustain baseline demand and buffer intermittent sources \cite{Jiang2005IPSN,Park2006SECON}. For the short, high-power bursts required by the laser-based transceiver, a supercapacitor delivers the peak current through a constant-current driver, thereby avoiding battery stress and voltage droop; this battery–supercapacitor hybrid is a standard solution for high-rate pulsed loads \cite{Shin2012JPS,Liu2022Photonics}. To coordinate these resources, the node uses an energy aware controller that supervises the battery and the supercapacitor.  As a result, the system follows a two stage cycle that alternates between energy accumulation and active communication, as outlined below.
\begin{itemize}
\item {Charging Phase:} The submarine illuminates the PV panel with a laser. The panel converts the received optical power into electrical energy that charges the lithium ion battery through a charge controller, and it can also precharge the supercapacitor when needed. The node remains in a low power state in which only housekeeping loads are active. This phase continues until the battery state of charge reaches a target that guarantees normal sensing tasks and an autonomous measurement window.
\item {Communication Phase:} After the energy targets are met, the node enables its optical receiver and the rest of the electronics. During laser based return link transmission the supercapacitor delivers the short duration burst power required by the transmitter, while the battery supports the baseline load. 
\end{itemize}

Although both phases are essential for autonomous operation, this work focuses on the communication phase with an emphasis on channel characterization and signal propagation. The analysis accounts for underwater optical impairments including absorption, scattering, turbulence induced fading, and pointing errors.

\begin{figure*}[b]
  \centering
\includegraphics[width=1\textwidth]{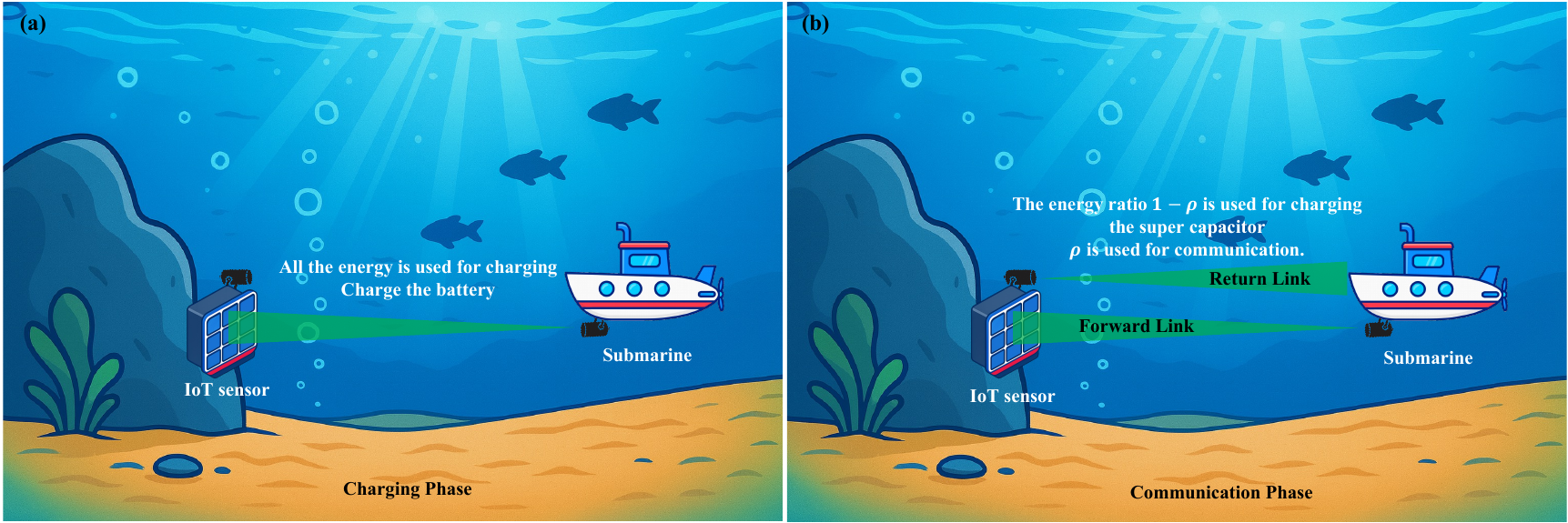}
  \caption{Diagram of the two phase underwater optical wireless communication system. (a) Charging phase.
 (b) Communication phase. }
  \label{fig:system_model}
\end{figure*}

Fig.~\ref{fig:system_model} illustrates the two phase operation of the proposed underwater optical wireless communication system. 

{ \subsection{Modeling Assumptions, Justification, and Applicability}
To keep the analysis tractable while capturing the dominant physics, the assumptions are listed with brief justification and applicability limits.

Assumption 1: PV panel as joint detector and harvester with a power split $\rho$. Optical SLIPT prototypes and circuit demonstrations show that a PV receiver can extract energy while conveying information by separating the direct current component for harvesting and the alternating component for detection \cite{Wang2015JSAC,Ammar2022CL}. This holds under small to moderate modulation depth and a stable operating point. Very strong ambient background or flicker may require additional optical filtering.

Assumption 2. Supercapacitor supplies the burst transmitter while the battery supports housekeeping loads. The battery and supercapacitor hybrid is common for pulsed loads because the supercapacitor delivers high peak current without stressing the battery and avoids voltage droop during short bursts \cite{Shin2012JPS,Liu2022Photonics}. This is valid when the burst duration is much smaller than the supercapacitor time constant and when the series resistance remains small over the used voltage window. Very long bursts need a detailed ladder model.

Assumption 3: Linearized models around the operating window for the supercapacitor and the PV front end. Supercapacitors are well represented by linear RC ladder equivalents within a chosen voltage span. When the PV junction is biased near a fixed operating point or under maximum power point tracking and the optical modulation depth is small, the photocurrent is approximately proportional to irradiance and the harvested power scales nearly linearly \cite{Tong2022OL,Lei2022JPHOT}. This holds for modest temperature drift and limited irradiance swing. Large temperature variation or large signal modulation requires nonlinear PV and thermal models.

Assumption 4: Independent forward and return links $h_1$ and $h_2$. The charging phase and the communication phase are separated in time by the supercapacitor and controller cycle. If this delay is longer than the channel coherence time the forward and return gains are statistically independent. 

Assumption 5. Piecewise constant bulk temperature and Gaussian jitter pointing with a finite aperture. The mixture EGG parameters summarize turbulence and bubble statistics measured under controlled thermal gradients \cite{Zedini2018EGG}. Residual pointing errors are captured by Gaussian beam jitter and an aperture coupling factor, which is widely used in optical propagation \cite{Farid2007Outage,Majumdar2015FSO}. The piecewise constant approximation holds over one communication cycle in typical shallow water tests. Large tilts, vignetting, or severe platform sway require refined geometry and tracking models.

Applicability across water conditions. The framework targets clear to coastal waters in the blue and green window where a line of sight path is feasible. For turbid harbor water or deep water with strong stratification the attenuation coefficients and mixture EGG parameters should be re identified to local conditions.}

\subsection{Signal and Power Models}
Let \(P_t\) denote the forward link optical power and $s_1(t)$ donate the signals from the submarine. So, the received signal  at the IoT sensor $y_s(t)$ can be modeled as
\begin{equation}
  y_s(t) = \rho\,\eta_r\,P_t\,h_1\,s_1(t) + n_1(t),
\end{equation}
and the instantaneous harvested  power $P_s(t)$ can be modeled as
\begin{equation}
  P_s(t) = (1-\rho)\,\eta_r\,P_t\,h_1,
\end{equation}
where \(n_1(t)\sim\mathcal{N}(0,\sigma_1^2)\), \(h_1\) is the forward link channel gain, and \(\eta_r\in(0,1)\) is the optical to electrical conversion efficiency of the PV front. The parameter \(\rho\in[0,1]\) is the fraction routed to the detector branch.

Then the node becomes the active transmitter and uses the instantaneous power \(P_s(t)\) to send return link data. The received return link signal at the submarine can be modeled as
\begin{align}
  y_u(t)
  &= \,\eta_s\,\eta_t\,P_s(t)\,h_2\,s_2(t) + n_2(t) \nonumber\\
  &= (1-\rho)\,\eta_t\,\eta_r\eta_s\,P_t\,h_1\,h_2\,s_2(t) + n_2(t),
\end{align}
where  \(\eta_t\) is the efficiency of converting electrical energy to light energy,  \(n_2(t)\sim\mathcal{N}(0,\sigma_2^2)\) and \(h_2\) is the return link channel gain, $\eta_s$ is the optical to electrical conversion efficiency for the receiver at the submarine. 

Because the charging and discharging of the supercapacitor introduce a time delay that typically exceeds the channel coherence time, \(h_1\) and \(h_2\) are modeled as statistically independent. Letting \(h = h_1 h_2\) denote the composite gain, the instantaneous return link SNR is
\begin{equation}
  \gamma_u = \overline{\gamma}_u\,h,
\end{equation}
with
\begin{equation}
  \overline{\gamma}_u = \frac{(1-\rho)\,\eta_t\,\eta_r\eta_s\,P_t}{\sigma_2^2}.
\end{equation}

\subsection{Channel Model}
\label{sec:channelmodel}
The end‐to‐end channel gain \(h_1\) from the submarine to the submerged sensor is influenced by three physical mechanisms along the underwater optical path. The first is the absorption and scattering loss, denoted by \(h_a\), which accounts for the exponential attenuation of optical power due to interactions with water molecules, particulate matter, and dissolved substances. The second component, \(h_t\), captures the irradiance fluctuations induced by dynamic underwater conditions such as air bubbles, turbulence, and temperature variations. Lastly, \(h_p\) represents the pointing error loss resulting from misalignment between the optical transmitter and receiver, which may be caused by surface wave motion or mechanical jitter. 

Assuming these effects are statistically independent, the overall channel gain is expressed as the product of the three components
\begin{equation}
  h_1 = h_a\, h_t\, h_p.
\end{equation}

For the reverse link, from the submerged sensor back to the submarine, the channel gain \(h_2\) is modeled using the same statistical framework. Under identical operating conditions (e.g., same wavelength, alignment characteristics, and water quality), \(h_1\) and \(h_2\) are assumed to be independent and identically distributed (i.i.d.). This i.i.d. assumption facilitates tractable system analysis while remaining consistent with practical deployment scenarios.

\subsubsection{Attenuation Loss}

Following the underwater attenuation model in \cite{Johnson2013AO}, the total attenuation of an optical signal propagating through seawater is primarily governed by two mechanisms: absorption and scattering. These effects are strongly dependent on the optical wavelength \(\lambda\), and their combined influence is quantified by the wavelength-dependent attenuation coefficient, which can be expressed as 
\begin{equation}
  c(\lambda) = a(\lambda) + b(\lambda),
\end{equation}
where \(a(\lambda)\) denotes the absorption coefficient, accounting for the conversion of optical energy into heat or chemical energy due to interactions with water molecules and dissolved organic matter, and \(b(\lambda)\) denotes the scattering coefficient, describing the redirection of photons due to suspended particles and microstructures in the water.

Based on Beer–Lambert's law, the total power attenuation over a propagation distance \(d_{vs}\) (from the submarine transmitter to the underwater sensor) is modeled as 
\begin{equation}
  h_a = \exp\left[-c(\lambda) d_{vs}\right],
\end{equation}
where $d_{vs}$ is the distance between the submarine and the sensor, \(h_a\) is the attenuation factor of the channel, representing the fraction of transmitted optical power that successfully reaches the receiver. The exponential form captures the cumulative effect of absorption and scattering along the propagation path.

To enable accurate analysis under different marine environments, typical values of the absorption and scattering coefficients \(a(\lambda)\) and \(b(\lambda)\) have been characterized for various water types using the Jerlov classification system. Table~\ref{Tab1} provides representative values of these coefficients at a commonly used wavelength of \(\lambda = 450\,\mathrm{nm}\), as adapted from \cite{Solonenko2015Jerlov}. This wavelength lies in the blue region of the visible spectrum, which is known to experience minimal attenuation in most oceanic waters, making it particularly favorable for underwater optical communication applications.
\begin{table}[ht]
\centering
\caption{Absorption and scattering coefficients for different water types at \(\lambda=450\) nm.}
\label{Tab1}
\scalebox{1}{
\begin{tabular}{l|cc}
\hline
Water Type                        & \(a(\lambda)\) (m\(^{-1}\)) & \(b(\lambda)\) (m\(^{-1}\)) \\
\hline
very clear ocean      & 0.008                               & 0.002                               \\
clear ocean          & 0.014                               & 0.003                               \\
intermediate clear    & 0.023                               & 0.004                               \\
coastal water         & 0.059                               & 0.009                               \\
turbid coastal water & 0.100                               & 0.020                               \\
\hline
\end{tabular}}
\end{table}

\begin{table}[!ht]
\centering
\caption{mixture EGG distribution parameters for different underwater conditions (Water Types 1–6).}
\label{Tab2_water_types}
\resizebox{\columnwidth}{!}{
\begin{tabular}{cccccccc}
\hline
\makecell{Water \\ Type} 
& \makecell{Temperature \\ Gradient \\ (\(^\circ\text{C}\cdot\text{cm}^{-1}\))} 
& \makecell{Bubbles \\ Level (L/min)} 
& \(\alpha\) 
& \(\beta\) 
& \(a\)    
& \(b\)    
& \(c\)      \\
\hline
1 & 0.050 & 2.4 & 0.21 & 0.33 & 1.4 & 1.2 & 17 \\
2 & 0.10  & 2.4 & 0.21 & 0.27 & 0.60 & 1.3 & 21 \\
3 & 0.15  & 2.4 & 0.18 & 0.16 & 0.23 & 1.4 & 23 \\
4 & 0.20  & 2.4 & 0.17 & 0.12 & 0.16 & 1.5 & 23 \\
5 & 0.050 & 4.7 & 0.46 & 0.34 & 1.0 & 1.6 & 36 \\
6 & 0.10  & 4.7 & 0.45 & 0.27 & 0.30 & 1.7 & 54 \\
\hline
\end{tabular}}
\end{table}

\subsubsection{Air‐Bubble and Temperature‐Induced Fading}

Small‐scale irradiance fluctuations in UWOC systems arise primarily due to the random scattering and refraction of light caused by microscopic air bubbles, particulate matter, and thermal inhomogeneities in the water. 

To accurately model the statistical behavior of these intensity fluctuations, Zedini \textit{et al.}~\cite{Zedini2018EGG} proposed the mixture EGG distribution, which captures both weak and strong turbulence regimes through a flexible mixture formulation. 
The PDF of the mixture EGG-distributed fading coefficient \(h_t\) is given by
\begin{equation}
  f_{h_t}(h_t)
  =  \frac{\alpha}{\beta}\exp\left(\frac{-h_t}{\beta}\right)
  + (1-\alpha) \frac{c h_t^{a c -1}}{b^{a c} \Gamma(a)}
    \exp\!\left[-\left(\frac{h_t}{b}\right)^c\right],
\label{PDFht}
\end{equation}
where \(\alpha \in [0,1]\) is the mixing coefficient that balances the contribution of the exponential and generalized Gamma components; \(\beta\) denotes the exponential decay parameter; and \(a\), \(b\), \(c\) are the shape and scale parameters of the generalized Gamma distribution. The exponential term captures the rapid fluctuations caused by scattering off bubbles and suspended particles, while the generalized Gamma term models the more persistent variations induced by thermal gradients and refractive index changes.

Table~\ref{Tab2_water_types} lists typical values of these parameters under different channel conditions (e.g., clear, coastal, or turbid water), as reported in~\cite{Zedini2018EGG}. This modeling approach enables a more accurate and tractable analysis of system performance metrics such as OP, bit error rate, and channel capacity under realistic underwater conditions.

\subsubsection{Pointing Error Loss}
The pointing error loss \(h_p\) quantifies the reduction in collected optical power when the beam centroid departs from the center of the receiver aperture because of mechanical vibrations, platform motion, or imperfect tracking.

We adopt the Gaussian jitter model \cite{Farid2007Outage,Majumdar2015FSO}. The lateral displacement of the beam spot at the receiver plane is modeled as bivariate Gaussian with zero mean and identical variance \(\sigma_s^2\) in both axes. The resulting radial displacement \(r\) follows a Rayleigh distribution
\begin{equation}
f_r(r)=\frac{r}{\sigma_s^2}\exp\!\left(-\frac{r^2}{2\sigma_s^2}\right), \quad r\ge 0 .
\end{equation}

Consider a circular aperture of radius \(r_a\) and a Gaussian beam with waist \(\omega_b\) at the receiver plane. The collected power for a given offset \(r\) is well approximated by \cite[Eq.~(9)]{Farid2007Outage}
\begin{equation}
h_p(r)\approx A_0 \exp\Bigl(-\frac{2r^2}{\omega_e^2}\Bigr),
\end{equation}
where $v=\frac{r_a\sqrt{\pi/2}}{\omega_b},\qquad
A_0=\operatorname{erf}^2(v),\qquad
\omega_e=\omega_b\sqrt{\frac{\sqrt{\pi}\,\operatorname{erf}(v)}{2\,v\,\exp(-v^2)}}$.

Applying a change of variables from \(r\) to \(h_p\) gives the probability density function of the pointing loss
\begin{equation}
f_{h_p}(h_p)=\frac{\mu_s^2}{A_0^{\mu_s^2}}\,h_p^{\mu_s^2-1},\qquad 0\le h_p\le A_0 ,
\label{PDFhp}
\end{equation}
with \(\mu_s=\omega_e/(2\sigma_s)\).Here \(A_0\) is the maximum collected power fraction under perfect alignment and \(\mu_s\) captures the ratio between the effective beam radius and the jitter level.

\subsubsection{Composite Channel Model}
Based on the statistical modeling of attenuation, turbulence-induced fading, and pointing error loss, we now derive the closed-form expression for the PDF of the forward link channel gain \(h_1 = h_a h_t h_p\). The PDF of \(h_1\) is summarized in the following theorem.
\begin{theorem}
\label{thm:pdf_h1}
The PDF of the  channel gain from the submarine to the submerged sensor  \(h_1 = h_a h_t h_p\) is given by
\begin{equation}
\label{PDFh1}
\begin{aligned}
f_{h_1}(h)& 
= \frac{\alpha \mu_s^2}{h}{\rm G}_{1,2}^{2,0}\left[\frac{h}{A_0 h_a \beta}\Bigm|\substack{\mu_s^2+1\\1,\mu_s^2}\right]
\\&+ \frac{(1-\alpha) \mu_s^2}{\Gamma(a) h} {\rm G}_{1,2}^{2,0}\left[\left(\frac{h}{A_0 h_a b}\right)^c\Bigm|\substack{\tfrac{\mu_s^2}{c}+1\\a, \tfrac{\mu_s^2}{c}}\right],
\end{aligned}
\end{equation}
where \(G_{p,q}^{ m,n}[\cdot]\) is the Meijer‐G function.
\begin{proof}
See Appendix~\ref{proof_pdf_h1}.
\end{proof}
{\noindent{Remarks:}
The Meijer G function is widely used in analytical modeling of FSO channels because it unifies many special functions and probability laws under a single framework \cite{shang2024enhancing,zedini2014performance,shang2025novel,shangOIRS}. It leads to compact closed form expressions for key performance metrics such as OP, average BER, and ergodic capacity. This representation is also practical for computation. Mathematica provides a built in Meijer-G implementation with arbitrary precision arithmetic.}
\end{theorem}

\section{Performance Analysis}
\label{Perana}
In this section, we analyze the performance of the proposed UWOC system by deriving key metrics such as the harvested power, SNR distribution, OP, average BER, and ergodic capacity. 

\subsection{Harvested Power Analysis}
To begin with, we focus on the energy harvesting aspect of the system, which is fundamental to ensuring autonomous and sustainable operation of the underwater sensor node.
\begin{corollary}
The PDF of \(P_s\) is given as
\begin{equation}
\label{PDFPs}
\begin{aligned}
f_{P_s}(p) &= \frac{\alpha \mu_s^2}{p} {\rm G}_{1,2}^{2,0}\left[\frac{p}{A_0 h_a \beta\, (1-\rho)\eta_r P_t}\Bigm|\substack{\mu_s^2+1\\1,\mu_s^2}\right] 
\\&+ \frac{(1-\alpha) \mu_s^2}{\Gamma(a)p} {\rm G}_{1,2}^{2,0}\left[\left(\frac{p}{A_0 h_a b\, (1-\rho)\eta_r P_t}\right)^c\Bigm|\substack{\tfrac{\mu_s^2}{c}+1\\a, \tfrac{\mu_s^2}{c}}\right],
\end{aligned}
\end{equation}
{ \begin{proof}
From the system definition $P_s=(1-\rho)\,\eta_r\,P_t\,h_1$. By change of variables,
\begin{equation}
f_{P_s}(p)=\frac{1}{(1-\rho)\eta_r P_t}\,
f_{h_1}\!\left(\frac{p}{(1-\rho)\eta_r P_t}\right),\qquad p>0.
\end{equation}
Substitute the expression of $f_{h_1}(\cdot)$ from~\eqref{PDFh1}, \eqref{PDFPs} can be obtained, which  completes the proof.
\end{proof}}
\end{corollary}

We now derive the average harvested power, which directly determines how long the node must remain in the charging phase before it can enter communication. A higher average harvested power shortens the time needed to precharge the supercapacitor and replenish battery headroom for housekeeping, increases the achievable duty cycle, and reduces access latency. It also exposes clear design tradeoffs through the split factor \(\rho\), the transmit power, the optical–electrical conversion efficiency, alignment quality, and water clarity, so it is a central quantity for scheduling, energy budgeting, and reliability in long duration deployments.
\begin{corollary}
The average harvested power at the sensor in the wake-up phase is given in closed form as
\begin{equation}
\label{E_Ps}
\mathbb{E}[P_s]
= \frac{A_0h_a(1-\rho)\eta_r{P}_t\mu_s^2}{\left(1 + \mu_s^2\right)}
\Biggl[
  \alpha\beta
  + \frac{(1-\alpha)b\Gamma\left(a + \tfrac{1}{c}\right)}
         {\Gamma(a)}
\Biggr].
\end{equation}
{\begin{proof}
 Start from the PDF in \eqref{PDFPs} and compute the first moment
\begin{equation}
\mathbb{E}[P_s]=\int_{0}^{\infty} p\, f_{P_s}(p)\,dp.
\end{equation}
Then, using \cite[Eq.~(2.8)]{Htran2}, which is given as
\begin{equation}
\label{eq:meijerG-mellin}
\begin{aligned}
&\int_{0}^{\infty} x^{s-1}\,
{\rm G}_{p,q}^{m,n}\!\left(
  \beta\, x^{c}\,\Biggm|\,
  \begin{matrix}
    a_1,\dots,a_p\\
    b_1,\dots,b_q
  \end{matrix}
\right)\,dx
\\&= \frac{1}{c}\,\beta^{-s/c}\,
\frac{
\prod_{j=1}^{m}\Gamma\!\left(b_j+\tfrac{s}{c}\right)\,
\prod_{i=1}^{n}\Gamma\!\left(1-a_i-\tfrac{s}{c}\right)
}{
\prod_{j=m+1}^{q}\Gamma\!\left(1-b_j-\tfrac{s}{c}\right)\,
\prod_{i=n+1}^{p}\Gamma\!\left(a_i+\tfrac{s}{c}\right)
},
\end{aligned}
\end{equation}
the expectation of $P_s$ can be obtained as \eqref{E_Ps}.
\end{proof}}
\end{corollary}

\subsection{PDF of the Composite Channel Gain}
Next, we turn our attention to the statistical modeling of the composite channel gain, which directly affects the received signal strength and SNR. The PDF of the composite channel gain \( h = h_1 h_2 \), which accounts for both turbulence-induced fading and pointing errors, plays a pivotal role in accurately characterizing the statistical behavior of received SNR. 
\begin{theorem}\label{thm:PDFh}
The PDF of the composite channel gain 
\( h = h_{1}h_{2}\)  is given by
\begin{equation}
\begin{aligned}
&f_{h}(h)
= \frac{(\alpha \mu_{s}^2)^{2}}{h} 
  {\rm G}_{2,4}^{4,0}\Biggl[
    \frac{h}{(A_{0}h_{a}\beta)^{2}}
    \Bigm|\substack{\mu_{s}^{2}+1, \mu_{s}^{2}+1\\
                     1, \mu_{s}^{2}, 1, \mu_{s}^{2}}
  \Biggr] \\&
+ \frac{2\alpha(1-\alpha) \mu_{s}^{4}}{\Gamma(a) c^{2} h} 
  {\rm H}_{2,4}^{4,0}\Biggl[
    \frac{h}{(A_{0}h_{a})^{2}\beta b}
    \Bigm|\substack{(\tfrac{\mu_{s}^{2}}{c}+1,\tfrac{1}{c})(\tfrac{\mu_{s}^{2}}{c}+1,\tfrac{1}{c})
                  \\(a,\tfrac{1}{c})(\tfrac{\mu_{s}^{2}}{c},\tfrac{1}{c}), (1,1), (\tfrac{\mu_{s}^{2}}{c},\tfrac{1}{c})}
  \Biggr] \\&
+ \left[\tfrac{(1-\alpha) \mu_{s}^{2}}{\Gamma(a)}\right]^{2}
  \frac{1}{c h} 
  {\rm G}_{2,4}^{4,0}\Biggl[
    \frac{h}{(A_{0}h_{a} b)^{2}}
    \Bigm|\substack{\tfrac{\mu_{s}^{2}}{c}+1, \tfrac{\mu_{s}^{2}}{c}+1\\
                     a, \tfrac{\mu_{s}^{2}}{c}, a, \tfrac{\mu_{s}^{2}}{c}}
  \Biggr].
\end{aligned}
\end{equation}
\label{PDFh}
\begin{proof}
See Appendix~\ref{app:proof_pdf_h}.
\end{proof}
\end{theorem}

\subsection{PDF and CDF of SNR}
Having established the composite fading statistics, we now derive the distribution of the instantaneous SNR \(\gamma_u\), which serves as the basis for computing other critical metrics such as OP and average BER. 
\begin{theorem}\label{thm:pdf_gamma_u}
Then the  PDF and CDF of \(\gamma_u\) are
\begin{equation}
\label{eq:pdf_gamma_u_full}
\begin{aligned}
f_{\gamma_u}&(\gamma)
= \frac{(\alpha\mu_s^2)^2}{\gamma} 
  {\rm G}_{2,4}^{4,0}\Biggl[
    \frac{\gamma{\overline{\gamma}}_u^{-1}}{ (A_0h_a\beta)^2}
    \Biggm|\substack{\mu_s^2+1, \mu_s^2+1\\
                     1, \mu_s^2, 1, \mu_s^2}
  \Biggr] \\&+
 2\alpha \frac{(1-\alpha) \mu_s^4}{\Gamma(a) c^2 \gamma}
  {\rm H}_{2,4}^{4,0}\Biggl[\frac{\gamma\left({\overline{\gamma}}_u\beta b\right)^{-1}}{ (A_0h_a)^2 }
    \Biggm|\substack{(\tfrac{\mu_{s}^{2}}{c}+1,\tfrac{1}{c})(\tfrac{\mu_{s}^{2}}{c}+1,\tfrac{1}{c})
                  \\(a,\tfrac{1}{c})(\tfrac{\mu_{s}^{2}}{c}\tfrac{1}{c})(1,1)(\tfrac{\mu_{s}^{2}}{c}\tfrac{1}{c})}
  \Biggr]\\&+ 
  \left[\tfrac{(1-\alpha) \mu_s^2}{\Gamma(a)}\right]^2
  \frac{1}{c \gamma} 
  {\rm G}_{2,4}^{4,0}\Biggl[\left[
    \frac{\gamma\overline{\gamma}_u^{-1}}{(A_0h_a b)^2}\right]^c
    \Biggm|\substack{\tfrac{\mu_s^2}{c}+1, \tfrac{\mu_s^2}{c}+1\\
                     a, \tfrac{\mu_s^2}{c}, a, \tfrac{\mu_s^2}{c}}
  \Biggr].
\end{aligned}
\end{equation}
\label{PDFgammau}
\begin{equation}
\begin{aligned}
&F_{\gamma_u}(\gamma)
= (\alpha\mu_s^2)^2
  {\rm G}_{3,5}^{4,1}\Biggl[
    \frac{\gamma{\overline{\gamma}}_u^{-1}}{ (A_0h_a\beta)^2}
    \Biggm|\substack{1,\mu_s^2+1, \mu_s^2+1\\
                     1, \mu_s^2, 1, \mu_s^2, 0}
  \Biggr]\\& + 
  \frac{2 \alpha(1-\alpha) \mu_s^4}{\Gamma(a) c^2 }
  {\rm H}_{3,5}^{4,1}\left[\frac{\gamma{\left({\overline{\gamma}}_u\beta b\right)}^{-1}}{ (A_0h_a)^2 }\!\!
    \Biggm|\!\!\substack{(1,1)
      (\tfrac{\mu_s^2}{c}+1,\tfrac1c) (\tfrac{\mu_s^2}{c}+1,\tfrac1c)\\(a,\tfrac1c)
      (\tfrac{\mu_s^2}{c},\tfrac1c)(1,1) (\tfrac{\mu_s^2}{c},\tfrac1c)(0,1)
    }
  \right]\\&+ 
  \left[\tfrac{(1-\alpha) \mu_s^2}{\Gamma(a)c}\right]^2
  {\rm G}_{3,5}^{4,1}\left[\left[
    \frac{\gamma\overline{\gamma}_u^{-1}}{(A_0h_a b)^2}\right]^c
    \Biggm|\substack{1,\tfrac{\mu_s^2}{c}+1, \tfrac{\mu_s^2}{c}+1\\
                     a, \tfrac{\mu_s^2}{c}, a, \tfrac{\mu_s^2}{c}, 0}\!\right].
\end{aligned}
\label{CDFgammau}
\end{equation}
\begin{proof}
Follows immediately by change of variables \(h=\frac{\gamma_u}{\overline{\gamma}_u}\) \eqref{PDFgammau} can be obtained. Starting from the definition of the cumulative distribution function,
\begin{equation}
  F_{\gamma_u}(\gamma_u)
  = \int_{0}^{\gamma_u} f_{\gamma_u}(x) \mathrm{d}x.
\label{eq:app_cdf_def}
\end{equation}
By substituting the closed‐form expression for \(f_{\gamma_u}(x)\) from \eqref{PDFgammau} into \eqref{eq:app_cdf_def}, and then applying the Mellin‐Barnes integral identity for Fox–H functions \cite[Eq.~2.25.2.2]{prudnikov1}, one directly arrives at the result stated in \eqref{CDFgammau}. This completes the proof.
\end{proof}
\end{theorem}

{ The diversity order captures the slope of the outage curve in the high SNR regime. It is insensitive to coding gain and other constant factors, so it provides a clean measure of how fast the OP decays as the average SNR $\overline{\gamma}_u$ increases. This helps compare the contributions of the components in \eqref{CDFgammau}. For any fixed threshold $\gamma>0$,
\begin{equation}
\label{diversity}
\mathcal{G}_{d}= -\lim_{\overline{\gamma}_u\to\infty}\frac{\log F_{\gamma_u}(\gamma)}{\log \overline{\gamma}_u}
=\min\left\{\,1,\ \mu_s^{2},\ ac\,\right\}.
\end{equation}

\begin{proof}
For $\overline{\gamma}_u\to\infty$ the arguments of the involved ${\rm G}$ and ${\rm H}$ functions go to zero. By the standard small argument expansions of Meijer–G and Fox–H functions \cite[Eq.~2.25.2.2]{prudnikov1}, each block in \eqref{CDFgammau} can be written as a finite sum of terms proportional to
\begin{equation}
\overline{\gamma}_u^{-d}\,[\ln \overline{\gamma}_u]^{m}\,(1+o(1)),
\end{equation}
where $m$ is a nonnegative integer that accounts for repeated poles, and the exponent $d$ belongs to the set $\{1,\mu_s^{2},ac\}$ obtained from the pole locations of the three blocks.where $m_i$ accounts for repeated poles and does not change the limiting slope. The possible exponents $d_i$ come from the pole locations as follows. The term with the smallest exponent dominates, while logarithmic factors are slower than any power and do not affect the limit, which proves the stated expression.
\end{proof} }


\subsection{Outage Probability}
Based on the derived CDF expression, we can now evaluate the OP, which quantifies the likelihood of communication failure due to insufficient SNR. The OP is defined as the probability that the instantaneous SNR at the receiver falls below a predefined threshold \(\gamma_{\mathrm{th}}\), which corresponds to the minimum SNR required to maintain an acceptable quality of service. In practical UWOC systems, this threshold may be determined by specific modulation and coding schemes or application-layer reliability requirements. A communication outage occurs whenever the received SNR dips below this critical threshold, resulting in data loss or decoding failure. Consequently, by evaluating the CDF in \eqref{CDFgammau} at \(\gamma=\gamma_{\mathrm{th}}\), i.e.,
\begin{equation}
    \mathrm{OP} = F_{\gamma_u}(\gamma_{\mathrm{th}}),
\end{equation}
It provides key insights into how different factors, such as water turbidity, pointing errors, and energy harvesting efficiency, impact the likelihood of successful data reception, which is essential for robust and sustainable sensor network design in underwater environments.

\subsection{SNR Moments}
We now derive the moments of the instantaneous SNR, which are essential for approximating average BER and ergodic capacity in scenarios where exact closed forms are difficult to obtain. Mathematically, the \( s \)-th moment of the instantaneous SNR, denoted as \( \mathbb{E}[\gamma_u^s] \), is defined as
\begin{equation}
\mathbb{E}\left[\gamma_u^s\right] = \int_0^{\infty} \gamma^s f_{\gamma_u}(\gamma) \, d\gamma,
\label{DefMoments}
\end{equation}
where \( f_{\gamma_u}(\gamma) \) is the PDF of \( \gamma_u \). 
\begin{corollary}
The $s$-th moments of $\gamma$ can be  expressed as
\begin{equation}
\begin{aligned}
\mathbb{E}[\gamma^{s}]
=\,& \left[\frac{\alpha    \mu_{s}^{2}    \Gamma(1+s)}{\mu_{s}^{2}+s}\right]^{2}
  \left[\frac{\overline\gamma_{u}^{-1}}{(A_{0}h_{a}\beta)^{2}}\right]^{s}\\&
+ \frac{2    \alpha    (1-\alpha)    \mu_{s}^{4}    \Gamma(a+s)    \Gamma(1+c    s)}
         {(\mu_{s}^{2}+c    s)^{2}    \Gamma(a)}
  \left[\frac{(\beta    b    \overline\gamma_{u})^{-1}}{(A_{0}h_{a})^{2}}\right]^{s}
\\&+ c
  \left[\frac{\Gamma(a+s)    (1-\alpha)    \mu_{s}^{2}}
              {\Gamma(a)    (\mu_{s}^{2}+c    s)}\right]^{2}
  \left[\frac{\overline\gamma_{u}^{-1}}{(A_{0}h_{a}b)^{2}}\right]^{c    s}.
\end{aligned}
\label{Momentsgammau}
\end{equation}
\begin{proof}

Substituing \eqref{PDFgammau} into \eqref{DefMoments} and using \cite[Eq.~(2.8)]{Htran2}, the moments of $\gamma$ can be evaluated according to (\ref{Momentsgammau}).
\end{proof}
\end{corollary}

\subsection{Average Bit-Error Rate}
With the SNR distribution and its moments ready, we now compute the average BER for M-PSK, M-QAM, and OOK. Although the link uses intensity modulation and direct detection (IM/DD), we adopt the standard subcarrier intensity modulation with direct detection (SIM-DD) model: the information symbols modulate an RF subcarrier electrically. The optical intensity carries this subcarrier \cite{Armstrong2009OFDMJLT,ArmstrongLowery2006PEO}. After square-law detection, the electrical passband is coherently demodulated. In SIM-DD, the post detection electrical SNR is proportional to the received optical power, so our instantaneous SNR model $\gamma_u=\overline{\gamma}_u\,h$ applies, with $\overline{\gamma}_u$ absorbing front end constants. A compact and unified expression of the average BER for various coherent M-QAM and M-PSK modulation schemes, as well as OOK modulation technique can be provided as 
\begin{equation}
\overline{P_{e}}=\delta_{B}\sum_{k=1}^{N_{B}}I\left(p_{B}, q_{B k}\right),
\end{equation}
where   $I\left(p_{B}, q_{B k}\right)$  is defined as 
\begin{equation}
\label{DefI}
\begin{aligned}
    I\left(p_{B},  q_{B k}\right)=&\frac{q_{Bk}^{p_{B}}}{2 \Gamma\left(p_{B}\right)} \int_{0}^{\infty} \gamma^{p_{B}-1} \exp \left(-q_{B k} \gamma\right) F_{\gamma}(\gamma) d \gamma,
\end{aligned}
\end{equation}
$N_{B}$, $\delta_{B}$ , $p_{B}$, and $q_{B k}$ are detailed in Table \ref{tab:my_label}.
\begin{table}[h]
\caption{Modulation Parameters}
\label{tab:my_label}
\centering
\footnotesize
\resizebox{\columnwidth}{!}{%
\begin{tabular}{lcccc}
\hline
Modulation & $\delta_{B}$ & $p_{B}$ & $q_{B_k}$ & $N_{B}$ \\ 
\hline
M-PSK   & $\displaystyle \frac{2}{\max(\log_2 M, 2)}$ 
        & $\tfrac12$ 
        & $\displaystyle \sin^2\!\left(\tfrac{(2k-1)\pi}{M}\right) \log_2 M$ 
        & $\displaystyle \max\left(\tfrac{M}{4}, 1\right)$ \\
M-QAM   & $\displaystyle \frac{4}{\log_2 M}\left(1-\frac{1}{\sqrt{M}}\right)$ 
        & $\tfrac12$ 
        & $\displaystyle \frac{3(2k-1)^2}{2(M-1)} \log_2 M$ 
        & $\displaystyle \frac{\sqrt{M}}{2}$ \\
OOK     & $1$ 
        & $\tfrac12$ 
        & $\tfrac12$ 
        & $1$ \\
\hline
\end{tabular}}
\end{table}

\begin{corollary}
\label{corollary32}
$I\left(p_{B}, q_{B k}\right)$ for our system is given as \eqref{Igammau}
\begin{figure*}[t]
\begin{equation}
\begin{aligned}
I\left(p_{B},q_{B_{k}}\right) \!
=\!\,&  \frac{(\alpha\mu_{s}^{2})^{2}}{2 \Gamma(p_{B})}
{\rm G}_{4,2}^{4,5}\Biggl[\frac{\left(\overline\gamma_{u} q_{B_{k}}\right)^{-1}}{(A_{0} h_{a} \beta)^{2}}
    \Bigg|\substack{1 - p_{B},1,\mu_{s}^{2}+1,\mu_{s}^{2}+1\\
                     1, \mu_{s}^{2}, 1, \mu_{s}^{2}, 0}
  \Biggr]\! +\!  \frac{\left[\tfrac{(1-\alpha) \mu_{s}^{2}}{\Gamma(a)c}\right]^{2}}{2 \Gamma(p_{B})} {\rm H}_{4,1}^{3,5}\Biggl[\left[
    \frac{\left(\overline\gamma_{u} q_{B_{k}}\right)^{-1}}{(A_{0} h_{a} b)^{2}}\right]^c \! 
    \Bigg|\substack{
      (1 - p_{B},c), (1,1), (\tfrac{\mu_{s}^{2}}{c}+1,1), (\tfrac{\mu_{s}^{2}}{c}+1,1)\\
      (a,1), (\tfrac{\mu_{s}^{2}}{c},1), (a,1), (\tfrac{\mu_{s}^{2}}{c},1), (0,1)
    }
  \Biggr] \\&+ \frac{\alpha (1-\alpha) \mu_{s}^{4} }
                {\Gamma(a) \Gamma(p_{B}) c^{2}}
  {\rm H}_{4,1}^{3,5}\Biggl[\frac{\left(\beta b \overline\gamma_{u} q_{B_{k}}\right)^{-1}}{(A_{0} h_{a})^{2}}
    \Bigg|\substack{
      (1 - p_{B},1), (1,1), (\tfrac{\mu_{s}^{2}}{c}+1,\tfrac1c), (\tfrac{\mu_{s}^{2}}{c}+1,\tfrac1c)\\
      (a,\tfrac1c), (\tfrac{\mu_{s}^{2}}{c},\tfrac1c), (1,1), (\tfrac{\mu_{s}^{2}}{c},\tfrac1c), (0,1)
    }
  \Biggr].
\end{aligned}
\label{Igammau}
\end{equation}
\hrule
\vspace{0.2cm}
\begin{equation}
\begin{aligned}
\overline{C}
= \,&        (\alpha\mu_s^2)^2
      {\rm G}_{3,5}^{5,1}\Biggl[
     \frac{\overline\gamma_u^{-1}}{(A_0   h_a   \beta)^2}
     \Bigg|\substack{
       0,  \mu_s^2+1,  \mu_s^2+1\\
       \mu_s^2,  1,  \mu_s^2,   0,   0
     }
   \Biggr]+ \frac{2   \alpha(1-\alpha)   \mu_s^4}{\Gamma(a)   c^2}
     {\rm H}_{3,5}^{5,1}\Biggl[
     \frac{(\beta   b   \overline\gamma_u)^{-1}}{(A_0   h_a)^2}
     \Bigg|\substack{
       (0,1),  (\tfrac{\mu_s^2}{c}+1,\tfrac1c),  (\tfrac{\mu_s^2}{c}+1,\tfrac1c)\\
       (a,\tfrac1c),  (\tfrac{\mu_s^2}{c},\tfrac1c),  (\tfrac{\mu_s^2}{c},\tfrac1c),  (0,1),  (0,1)
     }
   \Biggr]\\&
+ \left[\tfrac{(1-\alpha)   \mu_s^2}{\Gamma(a)}\right]^{2}\frac{1}{c^3}
     {\rm H}_{6,1}^{4,6}\Biggl[
   \frac{\overline\gamma_u^{-1}}{(A_0   h_a   b)^2   }
    \Bigg|\substack{
      (0,1),  (\tfrac{\mu_s^2}{c}+1,\tfrac1c),  (\tfrac{\mu_s^2}{c}+1,\tfrac1c),  (1,\tfrac1c)\\
      (a,\tfrac1c),  (a,\tfrac1c),  (\tfrac{\mu_s^2}{c},\tfrac1c),  (\tfrac{\mu_s^2}{c},\tfrac1c),  (0,\tfrac1c),  (0,1)
    }
  \Biggr].
\end{aligned}
\label{capacitygammau}
\tag{28}
\end{equation}
\hrule
\end{figure*}
\begin{proof}
See Appendix~\ref{app:proof_I_gammau}.
\end{proof}
\end{corollary}

\subsection{Ergodic Capacity}
Finally, we derive the ergodic capacity, a key indicator of the long-term average data rate achievable over the channel. Accordingly, the ergodic capacity of the end-to-end system is given by \cite{lapidoth}
\begin{align}
\label{DefC}
\overline{C}\triangleq \mathbb{E}[\ln(1+\gamma)]=\int_{0}^{\infty}\ln(1+ \gamma)f_\gamma(\gamma) d\gamma,
\end{align}
\begin{corollary}
\label{corollary33}
The ergodic capacity for the end-to-end system, $\overline{C}$, is given as \eqref{capacitygammau}.
\begin{proof}
See Appendix~\ref{app:proof_C_gammau}.
\end{proof}
\end{corollary}

\section{Numerical Results}
\label{sec:NUMERICAL ANALYSIS}

This section presents a detailed numerical evaluation of the proposed UWOC system to validate the analytical expressions and provide insights into its performance under practical underwater conditions.  Unless explicitly stated otherwise, all numerical results presented in this paper utilize these default parameters. 
\begin{table}[!ht]
\caption{System Parameters}
\begin{tabular}{cccc}
\hline
Parameters          & Values                      & Parameters     & Values      \\ \hline
$d_{vs}$               & 30 m                        & $\lambda$          & 450 nm    \\
$r_a$               & 5 cm                      & $\sigma_s$       & 0.5 $r_a$     \\
$\omega_b$               & 2 $r_a$                         & $a(\lambda)$          & 0.014 m$^{-1}$       \\
  $b(\lambda)$   & 0.003 m$^{-1}$  &           $\alpha$     & 0.213      \\
$\beta$           & 0.3291                         & $a$     & 1.4299        \\
$b$                   & 1.1817                       & $c$     & 17.1984     \\ $\eta_r$ &0.2 &$\eta_t$ & 0.8 \\ $\rho$ &0.8 &$\eta_s$ &0.9\\$\gamma_{\text{th}}$&2 dB& &  \\ \hline
\end{tabular}\label{tab2}
\end{table}

{ We simulate the UWOC link using Beer--Lambert attenuation with $a(\lambda)$ and $b(\lambda)$ at $\lambda=450\,\mathrm{nm}$ as in Table~\ref{Tab1}. This wavelength lies in the blue and green window where seawater loss is minimal, which is consistent with Jerlov measurements and ocean optics practice \cite{Johnson2013AO,Solonenko2015Jerlov}. Small scale irradiance fluctuations are generated with the mixture EGG model using the parameter sets in Table~\ref{Tab2_water_types}, which are fitted to controlled water tank data with calibrated temperature gradients and bubble injection \cite{Zedini2018EGG}. Pointing loss follows the Gaussian jitter model with a finite circular aperture. We set the aperture radius to $r_a=5\,\mathrm{cm}$ and the beam waist at the receiver plane to $\omega_b=2r_a$, then compute $A_0$, $\omega_e$, and $\mu_s=\omega_e/(2\sigma_s)$ from standard coupling formulas \cite{Farid2007Outage,Majumdar2015FSO}. The jitter standard deviation takes values $\sigma_s\in\{0.5,1,1.5,2\}\,r_a$ to cover mild and severe motion. Energy conversion efficiencies are fixed to $\eta_r=0.2$, $\eta_t=0.8$, and $\eta_s=0.9$. The split factor is set to $\rho=0.8$ to favor detection while preserving harvested power, and the distance set includes $d_{vs}\in\{10,30,50,70\}\,\mathrm{m}$ to represent short and mid range links in clear and coastal waters. Monte Carlo simulations use $N=10^{6}$ independent realizations with a fixed random seed. We generate $h_1=h_a h_t h_p$ and $h=h_1 h_2$ with independent hops, form $\gamma_u=\overline{\gamma}_u h$, and estimate PDFs, CDF, OP, average BER, and ergodic capacity. Analytical curves are evaluated with Meijer~G and Fox~H functions and are plotted together with the Monte Carlo results for validation.}

\subsection{Channel Characteristics}
\begin{figure}[!ht]
\centering\includegraphics[width=0.5\textwidth]{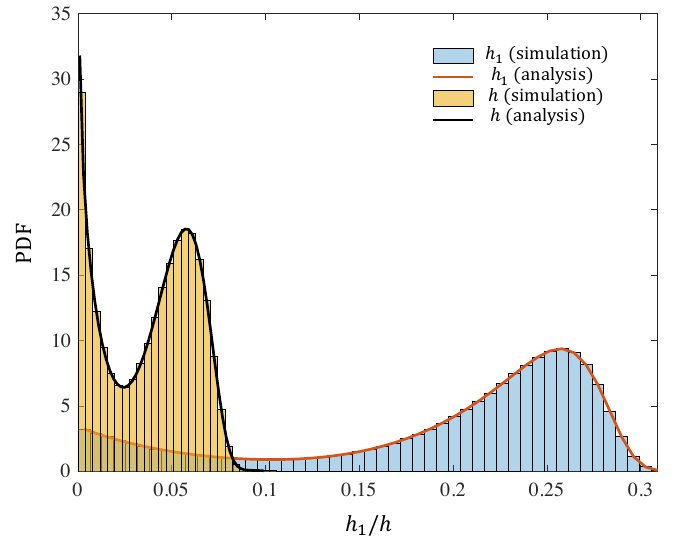}
\caption{Comparison between analytical and simulation results of the PDFs of single-link channel gain \(h_1\) and composite channel gain \(h = h_1 h_2\).}
  \label{figure2}
\end{figure}
This subsection investigates the statistical properties and energy transfer behavior of the underwater optical channel, focusing on both channel gain distribution and average harvested power under varying transmission conditions.

To validate the proposed analytical model, Fig.~\ref{figure2} compares the PDFs of the single-link channel gain \( h_1 \) and the composite channel gain \( h = h_1 h_2 \). The analytical results derived in this work are corroborated by MC simulations.  In addition to validation, Fig.~\ref{figure2} highlights the distinct statistical characteristics of \( h_1 \) and \( h \). Specifically, the PDF of the single-link gain \( h_1 \) follows a relatively unimodal distribution, whereas the composite channel gain \( h \) demonstrates a pronounced bimodal behavior. This difference stems from the multiplicative nature of independent random processes involved in the two-hop communication path, which introduces additional variance and structural changes in the composite gain distribution.

\begin{figure}[!ht]
  \centering
  \includegraphics[width=0.5\textwidth]{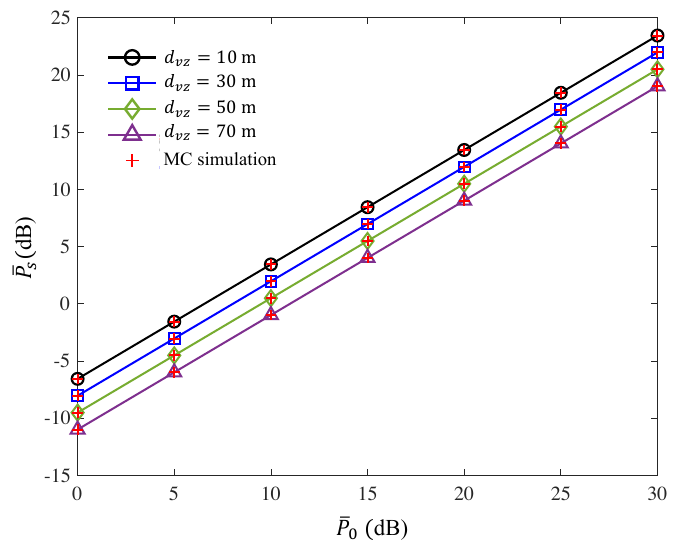}
\caption{Average harvested power at the underwater sensor \(\overline{P}_s\) versus average transmitted power from the submarine \(P_t\) for different link distances \(d_{vs}\).}
  \label{figure3}
\end{figure}

To better understand the energy transfer characteristics of the proposed UWOC system, we examine how the average harvested power at the sensor, denoted by \(\overline{P}_s\), varies with the average transmitted power from the submarine, \(P_t\), across different communication distances \(d_{vs}\), as illustrated in Fig.~\ref{figure3}. The results reveal an approximately linear growth of \(\overline{P}_s\) with increasing \(P_t\), which is consistent with the analytical expression provided in Corollary~2. As anticipated, larger transmission distances incur greater attenuation, resulting in reduced harvested power. For example, when \(P_t = 20\,\mathrm{dB}\), the corresponding values of \(\overline{P}_s\) are approximately \(13.4\,\mathrm{dB}\), \(12.0\,\mathrm{dB}\), \(10.5\,\mathrm{dB}\), and \(9.0\,\mathrm{dB}\) for \(d_{vs} = 10\,\mathrm{m}\), \(30\,\mathrm{m}\), \(50\,\mathrm{m}\), and \(70\,\mathrm{m}\), respectively. This clearly illustrates the adverse impact of propagation distance on the efficiency of energy harvesting.

\subsection{Performance metrics}
This subsection presents a comprehensive evaluation of the proposed UWOC system's performance, including OP, average BER, and ergodic capacity, under varying environmental and system parameters. Analytical results are compared against MC simulations to validate the accuracy of the proposed models.

\begin{figure}[!ht]
  \centering
  \includegraphics[width=0.5\textwidth]{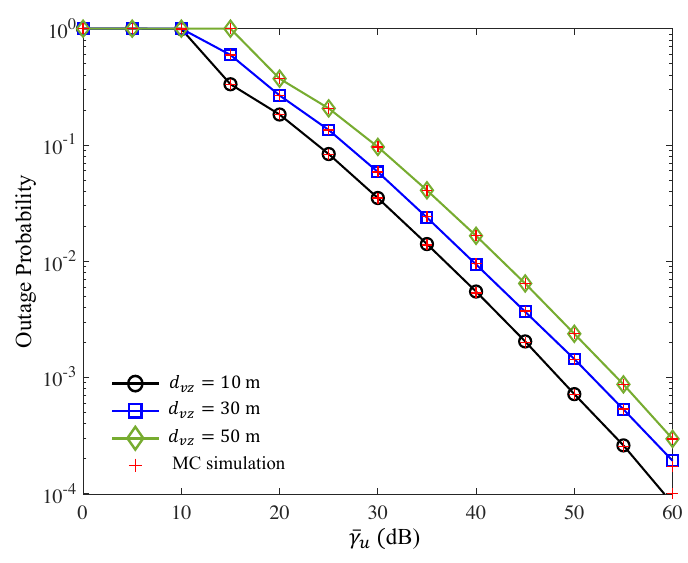}
\caption{OP of the return link transmission versus the average return link SNR \(\overline{\gamma}_u\) for different submarine–sensor distances \(d_{vs}\).}
  \label{figure4}
\end{figure}

To evaluate the reliability of the return link communication from the underwater sensor to the submarine, Fig.~\ref{figure4} illustrates the OP as a function of the average return link SNR, denoted by \(\overline{\gamma}_u\), under different values of the communication distance \(d_{vs}\). As observed, the OP decreases monotonically as \(\overline{\gamma}_u\) increases, indicating improved link reliability at higher SNRs. Moreover, shorter transmission distances consistently yield lower outage probabilities. For instance, when \(\overline{\gamma}_u = 30\,\mathrm{dB}\), the OP values are approximately \(3.5 \times 10^{-2}\), \(5.9 \times 10^{-2}\), and \(9.6 \times 10^{-2}\) for \(d_{vs} = 10\,\mathrm{m}\), \(30\,\mathrm{m}\), and \(50\,\mathrm{m}\), respectively. These findings highlight the significant impact of distance-dependent attenuation on the reliability of underwater optical return link transmission.

\begin{figure}[!ht]
  \centering
  \includegraphics[width=0.5\textwidth]{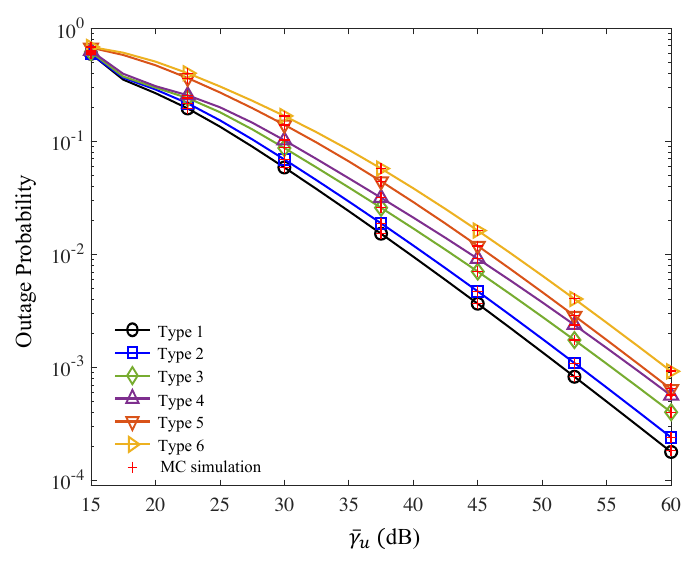}
    \caption{OP of the return link transmission versus the average return link SNR \(\overline{\gamma}_u\) under different underwater channel types (Type 1–6).}
\label{figure5}
\end{figure}

To further assess the robustness of the return link communication under varying underwater environmental conditions, we investigate the impact of different channel states characterized by water-induced turbulence and air bubbles. Fig.~\ref{figure5} presents the OP of the return link transmission as a function of the average return link SNR, \(\overline{\gamma}_u\), under six representative water types. These channel conditions, labeled as Water Types 1 through 6, reflect increasing levels of optical turbulence and scattering effects, governed by combinations of temperature gradient and bubble injection rate, as detailed in Table~\ref{Tab2_water_types}. Specifically, Water Types 1–4 are subject to a constant bubble level of 2.4~L/min, with increasing temperature gradients ranging from 0.05 to 0.20~\(^\circ\mathrm{C}\cdot\mathrm{cm}^{-1}\). In contrast, Water Types 5 and 6 introduce a higher bubble level of 4.7~L/min with temperature gradients of 0.05 and 0.10~\(^\circ\mathrm{C}\cdot\mathrm{cm}^{-1}\), respectively. As expected, the OP increases with more severe channel impairments. For example, at \(\overline{\gamma}_u = 30\,\mathrm{dB}\), the corresponding outage probabilities are approximately \(5.8 \times 10^{-2}\), \(6.8 \times 10^{-2}\), \(8.8 \times 10^{-2}\), \(1.0 \times 10^{-1}\), \(1.4 \times 10^{-1}\), and \(1.7 \times 10^{-1}\) for Water Types 1 through 6, respectively. This trend clearly confirms that both elevated temperature gradients and increased bubble levels intensify irradiance fluctuations, thereby degrading the reliability of underwater optical communication links.

\begin{figure}[!ht]
  \centering
  \includegraphics[width=0.5\textwidth]{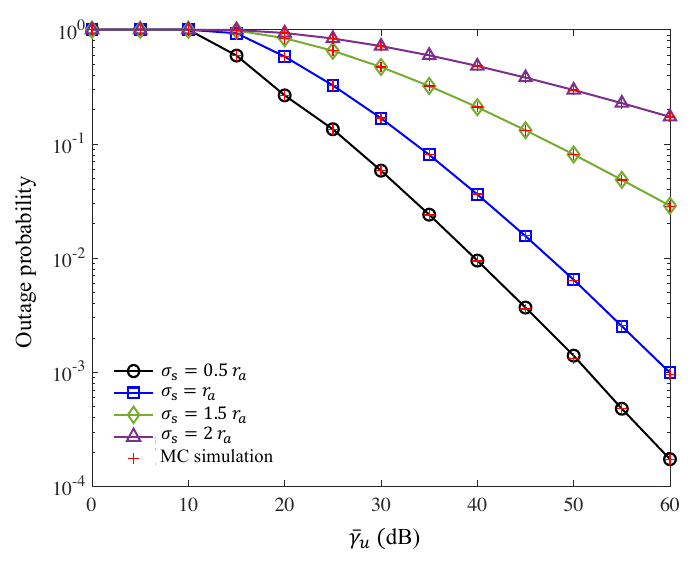}
\caption{OP versus average return link SNR \(\overline{\gamma}_u\) under different pointing error jitter levels.}
\label{figure6}
\end{figure}

To examine the influence of beam misalignment on the return link communication performance, we evaluate the impact of pointing error severity, characterized by the jitter standard deviation \(\sigma_s\), on the OP. Fig.~\ref{figure6} illustrates how the OP varies with the average return link SNR, \(\overline{\gamma}_u\), under different levels of pointing jitter. As depicted, smaller jitter values (e.g., \(\sigma_s = 0.5r_a\)) result in significantly enhanced link reliability, yielding notably lower outage probabilities. In contrast, larger jitter magnitudes (e.g., \(\sigma_s = 2r_a\)) cause substantial performance degradation due to increased beam misalignment. For instance, at \(\overline{\gamma}_u = 35\,\mathrm{dB}\), the corresponding OP values are approximately \(9.4\times10^{-3}\), \(3.6\times10^{-2}\), \(2.1\times10^{-1}\), and \(4.8\times10^{-1}\) for \(\sigma_s = 0.5r_a\), \(\sigma_s = r_a\), \(\sigma_s = 1.5r_a\), and \(\sigma_s = 2r_a\), respectively. These results highlight the critical importance of precise beam alignment in underwater optical wireless communication and reveal the detrimental impact that pointing errors can have on link reliability, especially in high-SNR regimes. 

\begin{figure}[!ht]
  \centering
  \includegraphics[width=0.5\textwidth]{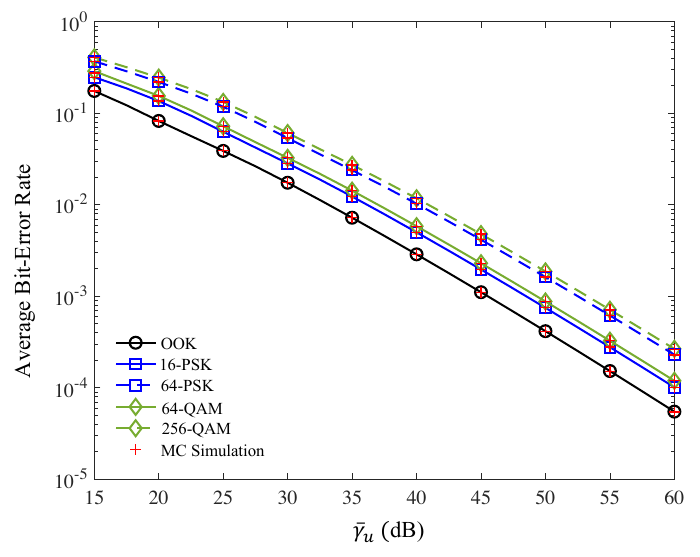}
\caption{Average BER versus average return link SNR \(\overline{\gamma}_u\) under various modulation schemes.}
\label{figure7}
\end{figure}

To investigate the effect of modulation format on return link reliability, we analyze the average BER performance under various modulation schemes. Fig.~\ref{figure7} shows the average BER as a function of the average return link SNR, \(\overline{\gamma}_u\), for five representative modulation formats: OOK, 16-PSK, 64-PSK, 64-QAM, and 256-QAM. As expected, the average BER decreases with increasing \(\overline{\gamma}_u\), and simpler modulation schemes demonstrate superior error performance due to their inherent resilience against underwater fading and noise. At \(\overline{\gamma}_u = 35\,\mathrm{dB}\), the average BER values for each scheme are approximately \(7.2 \times 10^{-3}\) for OOK, \(1.2 \times 10^{-2}\) for 16-PSK, \(1.4 \times 10^{-2}\) for 64-PSK, \(2.4 \times 10^{-2}\) for 64-QAM, and \(5.8 \times 10^{-2}\) for 256-QAM.  Overall, these results underscore the trade-off between spectral efficiency and robustness in the choice of modulation schemes for underwater optical return link transmission. 

\begin{figure}[!ht]
  \centering
  \includegraphics[width=0.5\textwidth]{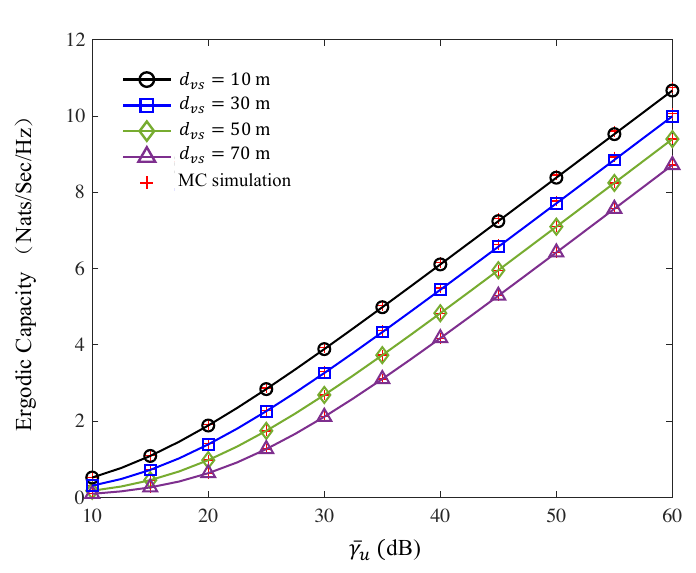}
\caption{Ergodic capacity versus average return link SNR \(\overline{\gamma}_u\) under different horizontal  distances \(d_{vs}\).}
\label{figure8}
\end{figure}

To assess the achievable data rates of the return link under varying propagation conditions, we examine the ergodic capacity as a function of the average return link SNR, \(\overline{\gamma}_u\), for different  distances. Fig.~\ref{figure8} presents the return link ergodic capacity under four representative values of \(d_{vs}\). As expected, the capacity increases monotonically with \(\overline{\gamma}_u\), reflecting the enhanced data throughput achievable with stronger signal power. Moreover, shorter horizontal  distances consistently yield higher capacity, owing to lower path loss and improved channel quality. For instance, at \(\overline{\gamma}_u = 40\,\mathrm{dB}\), the corresponding ergodic capacities are approximately 6.1, 5.4, 4.8, and 4.2 nats/sec/Hz for \(d_{vs} = 10\,\mathrm{m}\), \(30\,\mathrm{m}\), \(50\,\mathrm{m}\), and \(70\,\mathrm{m}\), respectively. These results emphasize the critical role of link distance in determining the information-theoretic limits of underwater optical wireless communication.

\begin{figure}[!ht]
  \centering
  \includegraphics[width=0.5\textwidth]{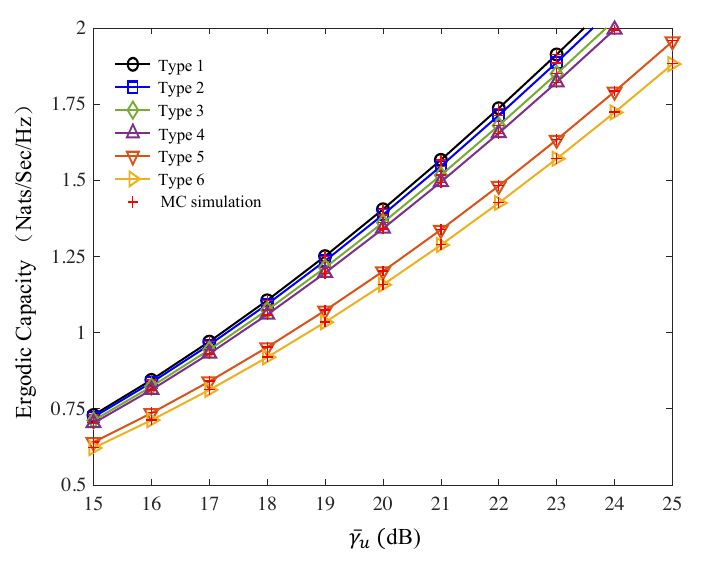}
 \caption{Ergodic capacity versus average return link SNR \(\overline{\gamma}_u\) under different water types.}
\label{figure9}
\end{figure}

To further evaluate the environmental sensitivity of the underwater optical wireless channel, we analyze the impact of different water types on the ergodic capacity of the return link transmission. Fig.~\ref{figure9} illustrates the ergodic capacity as a function of the average SNR, \(\overline{\gamma}_u\), across six representative water conditions characterized by varying temperature gradients and bubble injection rates, as defined in Table~\ref{Tab2_water_types}. As shown, the ergodic capacity increases with \(\overline{\gamma}_u\) for all water types, but the rate of increase and the absolute capacity values differ significantly depending on the water conditions. For example, at \(\overline{\gamma}_u = 20\,\mathrm{dB}\), the capacities are approximately 1.40, 1.38, 1.36, 1.34, 1.20, and 1.16 nats/sec/Hz for Water Types 1 through 6, respectively. This variation reflects the influence of both thermal and bubble-induced turbulence: Water Types 1–4 share a lower bubble level of 2.4~L/min, while Types 5 and 6 introduce a higher bubble level of 4.7~L/min, resulting in stronger scattering and absorption. These findings underscore the pronounced impact of environmental parameters on the capacity performance of underwater optical links.

\section{Conclusion and Future Research Direction}
\label{sec:Conclusion}
{ This paper developed an analytical framework for a two phase underwater optical link that realizes simultaneous lightwave information and power transfer with a photovoltaic panel receiver. We modeled distance dependent attenuation, turbulence induced fading with a mixture exponential–generalized Gamma law, and pointing loss under Gaussian jitter. Using Meijer G and Fox H functions, we derived closed form expressions for the channel gain and the instantaneous SNR together with OP, average bit error rate, ergodic capacity, and harvested power. Monte Carlo simulations verified the analysis and quantified the effects of water type, pointing jitter, and link distance.

The results support design and operation. 
(i) The closed forms enable link budgeting that sets transmit power, wavelength, and distance to meet a target outage level under a specified water type. 
(ii) Receiver choices that increase the collection area and reduce jitter improve both reliability and energy intake through the parameters $A_0$ and $\mu_s$. 
(iii) SLIPT control can tune the split factor $\rho$ to trade detection SNR and harvested power, which shortens charging time and stabilizes the duty cycle. 
(iv) The formulas allow fast sensitivity studies without heavy simulations, which helps select robust operating points for field deployment.

 We will pursue three directions. 
(i) Machine learning based adaptive SLIPT control that learns the split factor and rate under time varying channels, with safety constraints informed by the analytical bounds. 
(ii) Multi node energy and information optimization under dynamic water conditions, including joint scheduling of charging and communication, cross layer rate allocation, and fairness across sensors. 
(iii) Experimental validation on underwater testbeds and on AUV platforms, including calibration of the mixture EGG and pointing parameters, evaluation of closed loop alignment, and integration with practical photovoltaic front ends.}

\section*{Acknowledgment}
The authors would like to express their sincere gratitude to Dr. Abla Kammoun from KAUST for valuable discussions. We also thank Dr. Yaoyao Zhang from KAUST for insightful comments that improved this work.

\appendices
\numberwithin{equation}{section}%

\section{Proof of Theorem~\ref{thm:pdf_h1}}
\label{proof_pdf_h1}
To derive the PDF of the  channel gain from the submarine to the underwater sensor, \(h_1 = h_a h_t h_p\), we first condition on the turbulence-induced fading \(h_t\).  Applying the change-of-variables formula gives  
\begin{equation}
  f_{h}(h)
  = \int_{h/(A_0 h_a)}^{\infty}
    f_{h_p}\left(\frac{h}{h_a h_t}\right) 
    \frac{f_{h_t}(h_t)}{h_a h_t} 
    \mathrm{d}h_t,
\label{eq:proof_pdf_h1_1}
\end{equation}
where \(f_{h_p}(\cdot)\) and \(f_{h_t}(\cdot)\) are, respectively, the PDFs of the pointing-error component and the turbulence-induced fading component. Substituting \eqref{PDFhp} and \eqref{PDFht} into \eqref{eq:proof_pdf_h1_1} yields  
\begin{equation}
\begin{aligned}
&f_h(h)
= \frac{\mu_s^2}{A_0^{\mu_s^2}}
   \frac{\alpha}{h_a  \beta}
  \left(\frac{h}{h_a}\right)^{\mu_s^2-1}
  \underbrace{\int_{\frac{h}{A_0 h_a}}^{\infty}
    h_t^{ \mu_s^2}
    \exp\!\left(-\frac{h_t}{\beta}\right)
     \mathrm{d}h_t}_{I_1}
\\
&+
 \frac{\mu_s^2  (1-\alpha)c}{A_0^{\mu_s^2} h_a b^{a c} \Gamma(a)}
  \left(\frac{h}{h_a}\right)^{ \mu_s^2 - 1}\!\!\!
  \underbrace{\int_{\frac{h}{A_0 h_a}}^{\infty}
    h_t^{a c - \mu_s^2 - 1}
    \exp\!\left[-\!\left(\tfrac{h_t}{b}\right)^c\right]
     \mathrm{d}h_t}_{I_2} .
\end{aligned}
\label{eq:proof_pdf_h1_2}
\end{equation}

To evaluate the integrals \(I_1\) and \(I_2\), we employ \cite[Eq.~(11)]{MeijerGalgorithm} together with the primary definition of the Meijer-\(G\) function \cite[Eq.~(9.301)]{intetable}.  { Specifically, recall that the elementary exponential admits a Meijer–$G$ representation \cite[Eq.~(11)]{MeijerGalgorithm}
\begin{equation}
\exp{(-x)}
= {\rm G}^{1,0}_{0,1}\!\left[x \,\Bigg|\; \begin{matrix}-\\[1pt] 0\end{matrix}\right],
\end{equation}
and by the primary Mellin–Barnes definition of the Meijer–$G$ function \cite[Eq.~(9.301)]{intetable} this can be written as the contour integral
\begin{equation}
\exp{(-x)}
= \frac{1}{2\pi i}\int_{\mathcal L}\Gamma(-s)\,x^{\,s}\,{\rm d}s,
\qquad x>0 .
\label{eq:exp_MB_basic}
\end{equation}}

Therefore, for \(I_1\), we rewrite the exponential term via its Mellin–Barnes representation:
\begin{equation}
I_1
= \frac{1}{2\pi i}
\int_{\mathcal{L}}
\Gamma(-s)\,\left(\tfrac{1}{\beta}\right)^{s}
\int_{\frac{h}{A_0 h_a}}^{\infty}
h_t^{s-\mu_s^2} \mathrm{d}h_t
 \mathrm{d}s,
\label{eq:proof_pdf_h1_3}
\end{equation}
and, in an analogous manner,  
\begin{equation}
I_2
= \frac{1}{2\pi i}
\int_{\mathcal{L}}
\Gamma(-s)\,\left(\tfrac{1}{b}\right)^{c s}
\int_{\frac{h}{A_0 h_a}}^{\infty}
h_t^{c s + a c - \mu_s^2 -1} \mathrm{d}h_t
 \mathrm{d}s.
\label{eq:proof_pdf_h1_4}
\end{equation}

Carrying out the inner integrations and simplifying, we obtain the closed-form Meijer-\(G\) expressions  
\cite[Eq.~(11)]{MeijerGalgorithm}:
\begin{equation}
I_1
= \frac{\alpha \mu_s^2}{A_0 h_a \beta}
\left(\tfrac{h}{h_a}\right)^{\mu_s^2-1}
{\rm G}_{0,1}^{1,0}\Biggl[
\tfrac{h}{A_0 h_a \beta}
\Bigm|\substack{\mu_s^2+1\\0, \mu_s^2}
\Biggr],
\label{eq:proof_pdf_h1_5}
\end{equation}
\begin{equation}
I_2
= \frac{1}{c}
\left(\tfrac{h}{A_0 h_a}\right)^{a c - \mu_s^2}
{\rm G}_{1,2}^{2,0}\Biggl[
\left(\tfrac{h}{A_0 h_a b}\right)^c
\Bigm|\substack{\tfrac{\mu_s^2}{c}-a+1\\0, \tfrac{\mu_s^2}{c}-a}
\Biggr].
\label{eq:proof_pdf_h1_6}
\end{equation}

Finally, substituting \eqref{eq:proof_pdf_h1_5} and \eqref{eq:proof_pdf_h1_6} back into \eqref{eq:proof_pdf_h1_2} and invoking the reduction formulas in \cite[Eq.~(1.59)]{Htran2} and \cite[Eq.~(1.60)]{Htran2}, we arrive at the closed-form expression given in \eqref{PDFh1}, which completes the proof.

\section{Proof of Theorem~\ref{thm:PDFh}}
\label{app:proof_pdf_h}

The composite channel gain is defined as \(h = h_1 h_2\), where \(h_1\) and \(h_2\) are independent and identically distributed (i.i.d.).  Conditioning on \(h_1\) and applying the change–of–variables formula gives  
\begin{equation}
f_{h}(h)
= \int_{0}^{\infty}
  \frac{1}{h_{1}} 
  f_{h_{1}}(h_{1}) 
  f_{h_{2}}\!\left(\tfrac{h}{h_{1}}\right) 
  \mathrm{d}h_{1}.
\label{app:proof_pdf_h1}
\end{equation}

Because \(h_2\) has the same PDF as \(h_1\), substituting \eqref{PDFh1} into \eqref{app:proof_pdf_h1} yields  
\begin{equation}
\begin{aligned}
f_h&(h)
= \frac{(\alpha\mu_s^2)^2}{h}
  \int_{0}^{\infty}
    \frac{1}{h_1} 
    {\rm G}_{1,2}^{2,0}\left[
      \tfrac{h_1}{A_0h_a\beta}
      \Bigm|\substack{\mu_s^2 + 1\\1,\mu_s^2}
    \right]\\&\times
    {\rm G}_{1,2}^{2,0}\left[
      \tfrac{h}{A_0 h_a \beta h_1}
      \Bigm|\substack{\mu_s^2 + 1\\1,\mu_s^2}
    \right]
   \mathrm{d}h_1
+ \frac{2\alpha (1-\alpha)\mu_s^4}{\Gamma(a) h}
  \int_{0}^{\infty}
    \frac{1}{h_1}\\&\times
    {\rm G}_{1,2}^{2,0}\left[
      \tfrac{h_1}{A_0h_a\beta}
      \Bigm|\substack{\mu_s^2 + 1\\1, \mu_s^2}
    \right]
    {\rm G}_{1,2}^{2,0}\!\left[
      \left(\tfrac{h}{A_0 h_a b h_1}\right)^{c}
      \Bigm|\substack{\tfrac{\mu_s^2}{c} + 1\\a,\tfrac{\mu_s^2}{c}}
    \right]
   \mathrm{d}h_1 \\&
+ \left[\tfrac{(1-\alpha) \mu_s^2}{\Gamma(a)}\right]^{2}
  \frac{1}{h}
  \int_{0}^{\infty}
    \frac{1}{h_1}
    {\rm G}_{1,2}^{2,0}\!\left[
      \left(\tfrac{h_1}{A_0h_ab}\right)^{c}
      \Bigm|\substack{\tfrac{\mu_s^2}{c} + 1\\a,\tfrac{\mu_s^2}{c}}
    \right]\\&\times
    {\rm G}_{1,2}^{2,0}\!\left[
      \left(\tfrac{h}{A_0h_abh_1}\right)^{c}
      \Bigm|\substack{\tfrac{\mu_s^2}{c} + 1\\a,\tfrac{\mu_s^2}{c}}
    \right]
   \mathrm{d}h_1.
\end{aligned}
\label{app:proof_pdf_h2}
\end{equation}

Invoking the primary definition of the Meijer-\(G\) function \cite[Eq.~(9.301)]{intetable} and expressing each Meijer-\(G\) function in its Mellin–Barnes form, we obtain  
\begin{equation}
\begin{aligned}
f_h(h)&    
= \frac{(\alpha\mu_s^2)^2}{h} \frac{1}{2\pi i}\int_{\mathcal{L}}
    \frac{\Gamma(1 - t) \Gamma(\mu_s^2 - t)}
         {\Gamma(1 + \mu_s^2 - t)}
    \left(\tfrac{h}{A_0 h_a b}\right)^{t}\\&\times\int_{0}^{\infty}
      h_1^{\,t-1}
      {\rm G}_{1,2}^{2,0}\left[
      \tfrac{h_1}{A_0h_a\beta}
      \Bigm|\substack{\mu_s^2 + 1\\1,\mu_s^2}
    \right]
    \mathrm{d}h_1
   \mathrm{d}t
\\&
+ \frac{2 \alpha \mu_s^4 (1-\alpha)}{\Gamma(a)h} \frac{1}{2\pi i}
  \int_{\mathcal{L}}
    \frac{\Gamma(a - t) \Gamma\!\left(\tfrac{\mu_s^2}{c}-t\right)}
         {\Gamma\!\left(1 + \tfrac{\mu_s^2}{c} - t\right)}
    \left(\tfrac{h}{A_0 h_a b}\right)^{ct}\\&\times
\int_{0}^{\infty}
      h_1^{\,ct-1}
     {\rm G}_{1,2}^{2,0}\left[
      \tfrac{h_1}{A_0h_a\beta}
      \Bigm|\substack{\mu_s^2 + 1\\1,\mu_s^2}
    \right]
    \mathrm{d}h_1
  \mathrm{d}t
\\&
+ \left[\tfrac{(1-\alpha) \mu_s^2}{\Gamma(a)}\right]^{2}
  \frac{1}{hc} \frac{1}{2\pi i}
  \int_{\mathcal{L}}
    \frac{\Gamma(a - t) \Gamma\!\left(\tfrac{\mu_s^2}{c}-t\right)}
         {\Gamma\!\left(1 + \tfrac{\mu_s^2}{c} - t\right)}
    \left(\tfrac{h}{A_0 h_a b}\right)^{ct}\\&\times
   \int_{0}^{\infty}
      h_1^{-ct -1} 
      {\rm G}_{1,2}^{2,0}\!\left[
      \left(\tfrac{h}{A_0h_abh_1}\right)^{c}
      \Bigm|\substack{\tfrac{\mu_s^2}{c} + 1\\a,\tfrac{\mu_s^2}{c}}
    \right]
    \mathrm{d}h_1
  \mathrm{d}t. 
\end{aligned}
\label{app:proof_pdf_h3}
\end{equation}

Applying the identity in \cite[Eq.~(2.8)]{Htran2} to each inner integral gives  
\begin{equation}
\begin{aligned}
&f_h(h)
=\frac{ (\alpha\mu_s^2)^2}{h2\pi i}
    \int_{\mathcal{L}}
      \frac{\Gamma(1 - t)\Gamma(\mu_s^2 - t)\Gamma(1 - t)\Gamma(\mu_s^2 - t)}
           {\Gamma(1 + \mu_s^2 - t)\Gamma(1 + \mu_s^2 - t)}\\&\times
      \left[\tfrac{h}{(A_0 h_a \beta)^2}\right]^{t}
    \mathrm{d}t + \frac{2 \alpha \mu_s^4(1-\alpha)}{ch2\pi i}
    \int_{\mathcal{L}}\left[\tfrac{h}{(A_0 h_a)^2\beta b}\right]^{c t}\Gamma(a - t)\\&\times
      \frac{\Gamma\!\left(\tfrac{\mu_s^2}{c}-t\right)\Gamma(1 - ct)\Gamma\!\left(\tfrac{\mu_s^2}{c}-t\right)}
           {\Gamma\!\left(1 + \tfrac{\mu_s^2}{c}-t\right)\Gamma\!\left(1 + \tfrac{\mu_s^2}{c}-t\right)}
    \mathrm{d}t + 
   \frac{\left[\tfrac{(1-\alpha) \mu_s^2}{\Gamma(a)}\right]^{2}}{c h2\pi i} \\
&\times
   \int_{\mathcal{L}}
      \frac{\Gamma(a - t) \Gamma\!\left(\tfrac{\mu_s^2}{c}-t\right)\Gamma(a - c t)\Gamma\!\left(\tfrac{\mu_s^2}{c}-t\right)}
           {\Gamma\!\left(1 + \tfrac{\mu_s^2}{c}-t\right)\Gamma\!\left(1 + \tfrac{\mu_s^2}{c}-t\right)}
      \left[\tfrac{h}{(A_0 h_a b)^2}\right]^{c t}
     \mathrm{d}t.
\end{aligned}
\label{app:proof_pdf_h4}
\end{equation}

Finally, invoking the reduction formula in \cite[Eq.~(1.59)]{Htran2} and simplifying each contour integral recovers the closed-form PDF given in \eqref{PDFh}, thus completing the proof.

\section{Proof for Corollary \ref{corollary32}}
\label{app:proof_I_gammau}

Using the Meijer-G function's definition in \cite[Eq.~(9.301)]{intetable}, \eqref{CDFgammau} can be written as
\begin{equation}
\begin{aligned}
F_{\gamma_u}&(\gamma)  
=  \frac{(\alpha\mu_s^2)^2}{2\pi i}\!
  \int_{\mathcal{L}}
    \frac{\Gamma(1 - t) \Gamma(\mu_s^2 - t) \Gamma(1 - t) \Gamma(\mu_s^2 - t) \Gamma(t)}
         {\Gamma(1 + \mu_s^2 - t) \Gamma(1 + \mu_s^2 - t) \Gamma(1 + t)}\\&\times
    \left[\frac{\gamma \overline\gamma_u^{-1}}{(A_0 h_a \beta)^2}\right]^{t}
   \mathrm{d}t+ \frac{2 \alpha (1-\alpha) \mu_s^4}{\Gamma(a) c^22\pi i}
  \int_{\mathcal{L}}\left[\frac{\gamma (\beta b \overline\gamma_u)^{-1}}{(A_0 h_a)^2}\right]^{t}
    \\&\times\frac{\Gamma(a - \tfrac{t}{c}) \Gamma\left(\tfrac{\mu_s^2}{c} -  \tfrac{t}{c}\right) \Gamma(1 -  t) \Gamma\!\left(\tfrac{\mu_s^2}{c} -  \tfrac{t}{c}\right) \Gamma(t)}
         {\Gamma\!\left(1 + \tfrac{\mu_s^2}{c} -  \tfrac{t}{c}\right) \Gamma\!\left(1 + \tfrac{\mu_s^2}{c} -  \tfrac{t}{c}\right) \Gamma(1 + t)}
   \mathrm{d}t\\&+ \left[\tfrac{(1-\alpha) \mu_s^2}{\Gamma(a)}\right]^{2}
   \frac{1}{c^2}
   \frac{1}{2\pi i}
  \int_{\mathcal{L}}\left[\frac{\gamma \overline\gamma_u^{-1}}{(A_0 h_a b)^2}\right]^{c t}\\&\times \frac{\Gamma(a - t) \Gamma\!\left(\tfrac{\mu_s^2}{c} - t\right) \Gamma(a - c t) \Gamma\!\left(\tfrac{\mu_s^2}{c} - t\right) \Gamma(t)}
         {\Gamma\!\left(1 + \tfrac{\mu_s^2}{c} - t\right) \Gamma\!\left(1 + \tfrac{\mu_s^2}{c} - t\right) \Gamma(1 + t)}
   \mathrm{d}t.
\end{aligned}
\label{app:proof_I_gammau1}
\end{equation}
Substituting \eqref{app:proof_I_gammau1}
 into \eqref{DefI}, we can obtain
\begin{equation}
\scalebox{0.9}{$\begin{aligned}
&I\!\left(p_{B},q_{B_{m}}\right)
= 
   \frac{q_{B_{m}}^{p_{B}}(\alpha\mu_{s}^{2})^{2}}{2 \Gamma(p_{B})2\pi i}
  \int_{\mathcal{L}}\int_{0}^{\infty}
      \gamma_{u}^{\, t + p_{B} - 1} 
      \exp\bigl(- q_{B_{m}}\gamma_{u}\bigr)\,\mathrm{d}\gamma_{u} \\&\times
      \frac{\Gamma(1 - t) \Gamma(\mu_{s}^{2} - t) \Gamma(1 - t) \Gamma(\mu_{s}^{2} - t) \Gamma(t)}
           {\Gamma\!\left(1 + \mu_{s}^{2} - t\right) \Gamma\!\left(1 + \mu_{s}^{2} - t\right) \Gamma(1 + t)}
    \left[\frac{\overline{\gamma}_{u}^{-1}}{(A_{0}h_{a}\beta)^{2}}\right]^{t}\mathrm{d}t\\&+\, \frac{\alpha (1-\alpha) \mu_{s}^{4}\,q_{B_{m}}^{p_{B}}}{\Gamma(a) c^{2} \Gamma(p_{B})}
  \frac{1}{2\pi i}
  \int_{\mathcal{L}} \int_{0}^{\infty}
      \gamma_{u}^{\, t + p_{B} -1} 
      \exp\bigl(- q_{B_{m}}\gamma_{u}\bigr)\,\mathrm{d}\gamma_{u}\\&\times
      \frac{\Gamma\!\bigl(a - \tfrac{t}{c}\bigr) \Gamma\!\bigl(\tfrac{\mu_s^2}{c} -  \tfrac{t}{c}\bigr) \Gamma(1 -  t) \Gamma\!\bigl(\tfrac{\mu_s^2}{c} -  \tfrac{t}{c}\bigr) \Gamma(t)}
           {\Gamma\!\bigl(1 + \tfrac{\mu_s^2}{c} -  \tfrac{t}{c}\bigr) \Gamma\!\bigl(1 + \tfrac{\mu_s^2}{c} -  \tfrac{t}{c}\bigr) \Gamma(1 + t)}
    \left[\frac{(\overline{\gamma}_{u}\beta b)^{-1}}{(A_{0}h_{a})^{2}}\right]^{ t}\mathrm{d}t
\\&+\, \left[\tfrac{(1-\alpha) \mu_{s}^{2}}{\Gamma(a)}\right]^{2}
    \frac{q_{B_{m}}^{p_{B}}}{c^{2} \Gamma(p_{B})}
  \frac{1}{2\pi i}
  \int_{\mathcal{L}}\int_{0}^{\infty}
      \gamma_{u}^{\, c t + p_{B}-1} 
      \exp\bigl(- q_{B_{m}}\gamma_{u}\bigr)\,\mathrm{d}\gamma_{u}\\&\times
      \frac{\Gamma(a - t) \Gamma\!\bigl(\tfrac{\mu_{s}^{2}}{c}-t\bigr) \Gamma(a - c t) \Gamma\!\bigl(\tfrac{\mu_{s}^{2}}{c}-t\bigr) \Gamma(t)}
           {\Gamma\!\bigl(1 + \tfrac{\mu_{s}^{2}}{c}-t\bigr) \Gamma\!\bigl(1 + \tfrac{\mu_{s}^{2}}{c}-t\bigr) \Gamma(1 + t)}
    \left[\frac{\overline{\gamma}_{u}^{-1}}{(A_{0}h_{a}b)^{2}}\right]^{c t}\mathrm{d}t. 
\label{app:proof_I_gammau2}
\end{aligned}$}
\end{equation}
Using the definition of Gamma Function $\Gamma(z)$ in \cite{intetable}, $\Gamma(z)=\int_0^{\infty} t^{z-1} e^{-t} d t$, \eqref{Igammau} can obtained.

\section{Proof for Corollary \ref{corollary33}}
\label{app:proof_C_gammau}
By employing the definition of the Meijer-G function from \cite[Eq.~(9.301)]{intetable}, the CDF given in \eqref{CDFgammau} can be rewritten as
\begin{equation}
\scalebox{0.9}{$\begin{aligned}
&f_{\gamma_u}(\gamma)
= \frac{(\alpha\mu_s^2)^2}{\gamma 2\pi i}  
  \int_{\mathcal{L}}
    \frac{\Gamma(1-t)  \Gamma(\mu_s^2 - t)  \Gamma(1-t)  \Gamma(\mu_s^2 - t)  \Gamma(t)}
         {\Gamma(1+\mu_s^2 - t)  \Gamma(1+\mu_s^2 - t)  \Gamma(1+t)}\\&\times
    \left[\frac{\gamma  \overline\gamma_u^{-1}}{(A_0  h_a  \beta)^2}\right]^{t}
    \mathrm{d}t 
+ \frac{2  \alpha  (1-\alpha)  \mu_s^4}{\Gamma(a)  c^2  \gamma2\pi i}   
  \int_{\mathcal{L}}\left[\frac{\gamma  (\beta  b  \overline\gamma_u)^{-1}}{(A_0  h_a)^2}\right]^{t}
    \frac{ \Gamma(t)}
         {\Gamma(1+t)}\\&\times
    \frac{\Gamma(a - t)  \Gamma\left(\tfrac{\mu_s^2}{c}-t\right)  \Gamma(1 - c  t)  \Gamma\left(\tfrac{\mu_s^2}{c}-t\right) }
         {\Gamma\left(1 + \tfrac{\mu_s^2}{c}-t\right)  \Gamma\!\left(1 + \tfrac{\mu_s^2}{c}-t\right)}\mathrm{d}t 
+   \frac{\left[\tfrac{(1-\alpha)  \mu_s^2}{\Gamma(a)}\right]^{2}}{c  \gamma2\pi i}   \int_{\mathcal{L}} \\&\times
    \frac{\Gamma(a - t)  \Gamma\left(\tfrac{\mu_s^2}{c}-t\right)  \Gamma(a - c  t)  \Gamma\!\left(\tfrac{\mu_s^2}{c}-t\right)  \Gamma(t)}
         {\Gamma\!\left(1 + \tfrac{\mu_s^2}{c}-t\right)  \Gamma\!\left(1 + \tfrac{\mu_s^2}{c}-t\right)  \Gamma(1+t)}
    \left[\frac{\gamma  \overline\gamma_u^{-1}}{(A_0  h_a  b)^2}\right]^{c t}\!\!
    \mathrm{d}t.
\end{aligned}$}
\label{app:proof_C_gammau1}
\end{equation}
Substituting \eqref{app:proof_I_gammau1} into \eqref{DefI}, we obtain the following integral representation for \(I(p_B, q_{B_m})\)
\begin{equation}
\begin{aligned}
&\overline{C}   
= \frac{(\alpha\mu_s^2)^2}{2\pi i}
   \int_{\mathcal{L}}
     \frac{\Gamma(1-t) \Gamma(\mu_s^2-t) \Gamma(1-t) \Gamma(\mu_s^2-t)}
          {\Gamma(1+\mu_s^2-t) \Gamma(1+\mu_s^2-t)}\\&\times
     \left[\frac{\overline\gamma_u^{-1}}{(A_0 h_a \beta)^2}\right]^{t}
     \int_{0}^{\infty}\gamma^{ t-1}\ln(1+\gamma) \mathrm{d}\gamma
    \mathrm{d}t \\&+ \frac{2 \alpha (1-\alpha) \mu_s^4}{\Gamma(a) c^2 2\pi i}
    \int_{\mathcal{L}}
      \frac{\Gamma(a-t) \Gamma\!\left(\tfrac{\mu_s^2}{c}-t\right) \Gamma(1-t) \Gamma\!\left(\tfrac{\mu_s^2}{c}-t\right)}
           {\Gamma\!\left(1+\tfrac{\mu_s^2}{c}-t\right) \Gamma\!\left(1+\tfrac{\mu_s^2}{c}-t\right)}\\&\times
      \left[\frac{(\beta b \overline\gamma_u)^{-1}}{(A_0 h_a)^2}\right]^{t}
      \int_{0}^{\infty}\gamma^{ t-1}\ln(1+\gamma) \mathrm{d}\gamma
     \mathrm{d}t \\&
  + \left[\tfrac{(1-\alpha) \mu_s^2}{\Gamma(a)}\right]^{2}
    \frac{1}{c 2\pi i}
    \int_{\mathcal{L}}
      \frac{\Gamma(a-t) \Gamma\!\left(\tfrac{\mu_s^2}{c}-t\right) \Gamma(a-t) \Gamma\!\left(\tfrac{\mu_s^2}{c}-t\right)}
           {\Gamma\!\left(1+\tfrac{\mu_s^2}{c}-t\right) \Gamma\!\left(1+\tfrac{\mu_s^2}{c}-t\right)}\\&\times
      \left[\frac{\overline\gamma_u^{-1}}{(A_0 h_a b)^2}\right]^{c t}
      \int_{0}^{\infty}\gamma^{ c t-1}\ln(1+\gamma) \mathrm{d}\gamma
     \mathrm{d}t.
\end{aligned}
\end{equation}
Finally, by applying the standard definition of the Gamma function, given by \(\gamma_u\) as the standard Gamma function
\(\Gamma(z) = \int_0^{\infty} t^{z - 1} e^{-t} \, \mathrm{d}t,\) as provided in \cite{intetable}, the closed-form result in \eqref{Igammau} is obtained, which completes the proof.

\bibliographystyle{IEEEtran}
\bibliography{IEEEexample}
\end{document}